\documentclass[12pt,leqno]{amsart}
\newtheorem{theo}{Theorem}[section]

\newcommand{\be}{\begin{equation}}
\newcommand{\ee}{\end{equation}}
\newcommand{\bea}{\begin{eqnarray}}
\newcommand{\eea}{\end{eqnarray}}
\newcommand{\beb}{\begin{eqnarray*}}
\newcommand{\eeb}{\end{eqnarray*}}
\usepackage{amssymb,amsfonts,amsthm,setspace,indentfirst}
\numberwithin{equation}{section}
\begin{document}
\title[Curvature Properties of Interior Black Hole Metric]{ {\normalsize \bf{Curvature Properties of Interior Black Hole Metric}}}
\author[R. Deszcz,  A.H. Hasmani, V.G. Khambholja, and A.A. Shaikh]{Ryszard Deszcz, Abdulvahid H. Hasmani, Vrajeshkumar G. Khambholja, and Absos Ali Shaikh$^*$}
\date{}
\address{\noindent\newline Ryszard Deszcz \newline Department of Mathematics\newline Wroc\l aw University of
Environmental and Life Sciences\newline Grunwaldzka 53, 50-357 Wroc\l aw, Poland}
\email{ryszard.deszcz@upwr.edu.pl}
\address{\noindent\newline Abdulvahid H. Hasmani \newline Department of Mathematics\newline Sardar Patel University 
\newline Vallabh Vidyanagar, India}
\email{ahhasmani@gmail.com}
\address{\noindent\newline Vrajeshkumar G. Khambholja \newline Department of Mathematics\newline B. V. M. Engineering College\newline Vallabh Vidyanagar, India}
\email{vrajassp@gmail.com}
\address{\noindent\newline Absos A. Shaikh \newline Department of Mathematics\newline University of
Burdwan, Golapbag\newline Burdwan-713104, West Bengal, India}
\email{aask2003@yahoo.co.in, \ \ aashaikh@math.buruniv.ac.in}
\begin{abstract}
A spacetime is a connected 4-dimensional semi-Riemannian manifold endowed with a metric tensor $g$ with signature $(- + + +)$. 
The geometry of a spacetime is described by the tensor $g$ and the Ricci tensor $S$ of type $(0, 2)$ 
whereas the energy momentum tensor of type $(0,2)$ describes the physical contents of the spacetime. Einstein's field equations relate $g$, $S$ and the energy momentum tensor and describe the geometry and physical contents of the spacetime. 
By solving Einstein's field equations for empty spacetime (i.e. $S = 0$) for a non-static spacetime metric, 
one can obtain the interior black hole solution, known as the interior black hole spacetime which infers that a remarkable change occurs in the nature of the spacetime, namely, the external spatial radial and temporal coordinates exchange their characters to temporal and spatial coordinates, respectively, and hence the interior black hole spacetime is a non-static one as the metric coefficients are time dependent. For the sake of mathematical generalizations, in the literature, there are many rigorous geometric structures constructed by imposing the restrictions to the curvature tensor of the space involving first order and second order covariant differentials of the curvature tensor. Hence a natural question arises  that which geometric structures are admitted by the interior black hole metric. The main aim of this paper is to provide the answer of this question so that the geometric structures admitting by such a metric can be interpreted physically.
\end{abstract}

\noindent\footnotetext{$^*$ Corresponding author\\
$\mathbf{2010}$\hspace{5pt}Mathematics\; Subject\; Classification: 53B20, 53B25, 53B50, 53C25, 53C40,  83C57.}
\maketitle
\flushbottom
\section{\bf Introduction}
In the theory of general relativity one of the exciting predictions is that there may exist regions of the spacetime, 
where the gravity is so strong that nothing not even light, can ever escape. Such regions are known as black hole of the spacetime. 
It is well known that the most general spherically symmetric, static, vacuum, asymptotically flat exact solution 
to Einstein's field equations is described by the Schwarzschild metric
\be\label{eq1.1}
ds^2 = -\left(1-\frac{2m}{z}\right)dt^2 + \left(1-\frac{2m}{z}\right)^{-1}dz^2 + z^2\left(d\theta^2 + \sin^2 \theta d\phi^2\right),
\ee
where $z$, $\theta$, $\phi$ are spherical polar coordinates, $t$ is the time measure by a clock at infinity and $m= \frac{MG}{c^2}$, $M$ being the mass of the central body, $G$ being the gravitational constant and $c$ being the velocity of light. 
Geometric properties of the Schwarzschild metric are presented among others in {\cite[Chapters 23-26] {Blau}},
{\cite[Section 8] {GrifPod}} and {\cite[Chapter 13] {B}}. 
We also refer to {\cite[Chapter 4.7] {Chen-2017}} for a short introduction to this subject. 
Curvature properties of pseudosymmetry type of this metric are given in \cite{{DeVerVra}, {SKS17}}.\\
\indent On the sphere $z = 2m$, the coefficient of $dz^2$ in (\ref{eq1.1}) tends to infinity and hence $z = 2m$ is a singularity. 
The metric has also a singularity at $z = 0$ as the coefficient at $dt^2$ in (\ref{eq1.1}) tends to infinity. We note that $z = 0$ is the centre of the spherical mass distribution of the star. Since the metric coefficients are coordinate dependent, $z = 0$ is a coordinate singularity, which implies that the inverse of the metric components $g_{33}$ and $g_{44}$ diverges even though at that point there is nothing such physical possibility. The geometry behaves very strangely when $z= 2m$ and it is called the Schwarzschild radius. The spherical surface associated with the Schwarzschild radius $z = 2m$ is null and corresponds to the black hole's ``event horizon'', where a freely falling particle can approach to the surface $z = 2m$ but never cross it and hence the inward falling particle need infinite time to reach the surface of the sphere $z = 2m$, which was first pointed out by Oppenheimer and Snyder \cite{OS}. We note that inside the Schwarzschild radius, $z$ and $t$ coordinates change their role in the sense that the $t$ coordinate becomes spacelike and $z$ coordinate becomes timelike. It is believed that the gravitational collapse of a compact body results in a singularity hidden beyond an event horizon. If the singularity were visible to the exterior region one would have a naked singularity which would open the realm for wild speculations \cite{ROSA}. This entails to Penrose's cosmic censorship conjecture \cite{PEN} which states that all physically reasonable spacetimes are globally hyperbolic, forbidding the existence of naked singularities, and only allowing singularities to be hidden behind event horizon. The realization that black hole could actually exists prompted a renewed interest in their mathematical properties and the last three decades have seen some remarkable developments in this respect. For details about the black hole in cosmology and astrophysics we refer the article of Carr \cite{CARR} and also references therein.\\
\indent The empty annular region of spacetime for a spherical star inside its Schwarzschild radius $2m$ and outside its physical radius 
$a$, $a < 2m$, that is for a black hole the spacetime geometry is characterized by the spherical symmetric non-static line element
\be\label{eq1.2}
ds^2 = -B(z, t)dt^2 + A(z, t)dz^2 + t^2\left(d\theta^2 + \sin^2 \theta d\phi^2\right),
\ee
where the coefficient functions can be obtained by solving Einstein's field equations for empty spacetime $S = 0$. The solution takes the form
\be\label{eq1.3}
ds^2 = -\left(\frac{2\xi}{t} - 1\right)^{-1}dt^2 + \left(\frac{2\xi}{t} - 1\right)dz^2 + t^2\left(d\theta^2 + \sin^2 \theta d\phi^2\right).
\ee
Thus (\ref{eq1.3}) represents metric for interior black hole 
\cite{{ROSA}, {HK}}, where $\xi$ may be determined from a direct confrontation with the exterior Schwarzschild solution. The metric (\ref{eq1.3}) is the interior black hole solution which represents the empty spacetime in the exterior region $z > a$ of a black hole and it is valid for $a$ such that $2m < z < a$. For physical significance and cosmological interpretation of the interior black hole solution we refer the reader \cite{ROSA} and also references therein. However in the interior black hole solution, a remarkable change occurs in the nature of spacetime namely the external spatial radial and temporal coordinates exchange their characters to temporal and spatial coordinates, respectively and hence the interior black hole solution is represented by a non-static spacetime as its metric coefficients are time dependent.\\
\indent The nature of a space is completely known by its curvature which can be explicitly determined by the metric of that space. 
In the literature of differential geometry there are several kinds of generalizations of various geometrical 
structures constructed by giving the curvature restrictions involving first and second order covariant derivatives.  
The notion of local symmetry is a generalization of the manifold of constant curvature and the study was initiated 
by Cartan in 1926 with full classification of such a space \cite{{CA1}, {CA2}}. 
Also a full classification of locally symmetric semi-Riemannian space is given by Cahen and Parker \cite{{CP1}, {CP2}}. 
During the last eight decades the process 
of generalization of locally symmetric spaces have been carried out by many authors around the globe in different directions, 
for instance, recurrent manifold by Walker \cite{WA}, 2-recurrent manifold by Lichnerowicz \cite{LI}, 
quasi-generalized recurrent manifold by Shaikh and Roy \cite{SR1}, hyper-generalized recurrent manifold 
by Shaikh and Patra \cite{SP}, weakly generalized recurrent manifold by Shaikh and Roy 
\cite{{SFR}, {SR2}}, semisymmetric manifold by Cartan \cite{CA}, 
pseudosymmetric manifold by Chaki \cite{CH}, pseudosymmetric manifold by Deszcz \cite{{DES}, {DEG}},
Ricci-pseudosymmetric manifold by Deszcz and Hotlo\'{s} \cite{DH}, 
weakly symmetric manifold by T$\acute{\mbox{a}}$massy and Binh \cite{TB}, weakly symmetric manifold by Selberg \cite{SE}.  
We mention that pseudosymmetry by Chaki and Deszcz are different and also weak symmetry of Selberg and T$\acute{\mbox{a}}$massy and Binh are different.\\
\indent We consider the semi-Riemannian manifold $(M,g)$, $\dim M = 4$,
equipped with the interior black hole metric given in (\ref{eq1.3}). 
Then this manifold is a interior black hole spacetime. 
The main subject of this paper is to investigate the geometric structures admitting by the interior black hole spacetime. 
The paper is organized as follows. In Section 2 we present 
definitions of some special tensors. In Section 3 we present basical facts on pseudosymmetric manifolds 
(in the sense of Deszcz \cite{LV1}). 
In Section 4 we deduce the curvature properties of interior black hole metric and found that 
interior black hole spacetime is a pseudosymmetric manifold.

Finally, in the last section (Appendix) we present the local components of the considered tensors of the metric 
(\ref{eq1.3}). We also mention that we have made all the calculations by a programme in Wolfram Mathematica. 
\section{\bf Some special tensors}
\indent Let $(M,g)$, $n = \dim M \ge 3$, be a connected smooth semi-Riemannian manifold 
with Levi-Civita connection $\nabla$ and semi-Riemannian metric $g$. 
For $(0,2)$-tensors $A$ and $B$ on $M$ 
we define their Kulkarni-Nomizu product $A\wedge B$ 
by (see, e.g., \cite{{D2}, {GLOG}, {SK14}, {SKgrt}, {SK16}, {SKA16}, {SRK16}})
\begin{eqnarray*}
(A\wedge B)(X_1,X_2,X,Y)&=&A(X_1,Y)B(X_2,X)+A(X_2,X)B(X_1,Y)\\
&-&A(X_1,X)B(X_2,Y)-A(X_2,Y)B(X_1,X).
\end{eqnarray*}
We define the endomorphisms $X\wedge_A Y$, $\mathcal{R}(X,Y)$,  
$\mathcal{C}(X,Y)$,  $\mathcal{P}(X,Y)$,  $\mathcal{W}(X,Y)$ and $\mathcal{K}(X,Y)$ 
by \cite{{D0}, {D11}, {D2}, {D9}, {D10}, {GLOG}, {GLOG1}}
\begin{eqnarray*}
&&(X\wedge_A Y)Z = A(Y,Z)X-A(X,Z)Y,\\
&&\mathcal{R}(X,Y)Z = [\nabla_X,\nabla_Y]Z-\nabla_{[X,Y]}Z,\\
&&\mathcal{C}(X,Y) = \mathcal{R}(X,Y)-\frac{1}{n-2}(X\wedge_g \mathcal{L} Y + \mathcal{L} X\wedge_g Y - \frac{r}{n-1}X\wedge_g Y),\\
&&\mathcal{P}(X,Y) = \mathcal{R}(X,Y)-\frac{1}{n-1}X\wedge_S Y,\\
&&\mathcal{W}(X,Y) = \mathcal{R}(X,Y)-\frac{r}{n(n-1)}X\wedge_g Y,\\
&&\mathcal{K}(X,Y) = \mathcal{R}(X,Y)-\frac{1}{n-2}(X\wedge_g \mathcal{L} Y + \mathcal{L} X\wedge_g Y),
\end{eqnarray*}
respectively, where $A$ is a $(0,2)$-tensor on $M$, $X,Y,Z\in\chi(M), \chi(M)$ being the Lie algebra of smooth vector fields on $M$. 
The Ricci tensor $S$, the Ricci operator ${\mathcal{L}}$ and the scalar curvature $r$ are defined by 
$S(X,Y) =  \mathrm{tr}\, \{ Z \mapsto {\mathcal{R}}(Z,X)Y \}$, 
$g({\mathcal L}X,Y) = S(X,Y)$ and 
$r  = \mathrm{tr}\, {\mathcal{L}}$, respectively.
We define the tensor $G$, the Riemannian-Christoffel curvature tensor $R$, the Weyl conformal curvature tensor $C$, 
the projective curvature tensor $P$, the concircular curvature tensor $W$ 
and the conharmonic curvature tensor $K$ of $(M,g)$, 
by
\cite{{D0}, {D11}, {D2}, {D9}, {D10}, {GLOG}, {GLOG1}}
\begin{eqnarray*}
&&
G(X_1, \ldots ,X_4) = g((X_1\wedge_g X_2)X_3,X_4) = \frac{1}{2}(g \wedge g)(X_1, \ldots ,X_4),\\
&&
R(X_1, \ldots ,X_4) = g(\mathcal{R}(X_1,X_2)X_3,X_4),\\
&&
C(X_1, \ldots ,X_4) = g(\mathcal{C}(X_1,X_2)X_3,X_4)\\
&&
= (R -\frac{1}{n-2}g \wedge S + \frac{r}{(n-1)(n-2)} G)(X_1, \ldots ,X_4),
\end{eqnarray*}
\begin{eqnarray*}
&&
P(X_1, \ldots ,X_4) = g(\mathcal{P}(X_1,X_2)X_3,X_4)\\
&&
=
R(X_1, \ldots ,X_4)
-
\frac{1}{n-1}( g(X_1,X_4) S(X_2,X_3) - g(X_2,X_4) S(X_1,X_3)),\\
&&
W(X_1, \ldots ,X_4) = g(\mathcal{W}(X_1,X_2)X_3,X_4) = (R - \frac{r}{n(n-1)}G)(X_1, \ldots ,X_4),\\
&&
K(X_1, \ldots ,X_4) = g(\mathcal{K}(X_1,X_2)X_3,X_4) = (R - \frac{1}{n-2} g \wedge S)(X_1, \ldots ,X_4)\\
&&
= (C - \frac{r}{(n-1)(n-2)} G) (X_1, \ldots ,X_4),
\end{eqnarray*}
respectively. 
For an $(0, k)$-tensor $T$, $k\geq 1$, and a symmetric $(0, 2)$-tensor $A$ we define the $(0, k+2)$-tensor $Q(A,T)$ 
\cite{{SK14}, {SKgrt}, {SKppsnw}, {SKppsn}} by
\begin{eqnarray*}
&&Q(A, T)(X_1,\ldots , X_k;X,Y)=((X\wedge_A Y) \cdot T)(X_1, \ldots , X_k)\\
&&=-T((X\wedge_A Y)X_1,X_2, \ldots , X_k)- \cdots - T(X_1, \ldots , X_{k-1},(X\wedge_A Y)X_k).
\end{eqnarray*}
The tensor $Q(A,T)$ is called the Tachibana tensor of the tensors $A$ and $T$, 
or shortly the Tachibana tensor \cite{{DGlogH}, {DGHSaw}}.
It is obvious that the tensor $Q(g,G)$ vanishes identically on any semi-Riemannian manifold. 
Therefore we have
$Q(g,R) = Q(g,W)$ and $Q(g,C) = Q(g,K)$. 
For an endomorphism $\mathcal{D}(X,Y)$ we define 
the $(0,4)$-tensor $D$ by
\begin{eqnarray*}
D(X_1, \ldots ,X_4) = g(\mathcal{D}(X_1,X_2)X_3,X_4) .
\end{eqnarray*}
Now for an $(0, k)$-tensor $T$, $k\geq 1$, and an endomorphism $\mathcal{D}(X,Y)$ we define 
the $(0, k+2)$-tensor $D \cdot T$ \cite{{SK14}, {SKgrt}, {SKppsnw}, {SKppsn}} by
\begin{eqnarray}
& & (D \cdot T) (X_1,\ldots, X_k;X,Y)= (\mathcal{D}(X,Y) \cdot T)(X_1, \ldots, X_k)\nonumber\\
&&=-T(\mathcal{D}(X,Y)X_1,X_2, \ldots, X_k)- \cdots - T(X_1, \ldots, X_{k-1}, \mathcal{D}(X,Y)X_k).
\label{ror}
\end{eqnarray}
Setting in the above formulas $\mathcal{D}(X,Y) = \mathcal{R}(X,Y),\ \mathcal{C}(X,Y),\ 
\mathcal{P}(X,Y),\ \mathcal{W}(X,Y),\ \mathcal{K}(X,Y)$, $T = R,\ S,\ C,\ P,\ W,\ K$ and $A = g$ or $S$, 
we obtain the tensors: 
$R\cdot R$, $R\cdot S$, $R\cdot C$, $R\cdot P$, $R\cdot W$, $R\cdot K$, 
$C\cdot R$, $C\cdot S$, $C\cdot C$, $C\cdot P$, $C\cdot W$, $C\cdot K$, 
$P\cdot R$, $P\cdot S$, $P\cdot C$, $P\cdot P$, $P\cdot W$, $P\cdot K$, 
$W\cdot R$, $W\cdot S$, $W\cdot C$, $W\cdot P$, $W\cdot W$, $W\cdot K$, 
$K\cdot R$, $K\cdot S$, $K\cdot C$, $K\cdot P$, $K\cdot W$, $K\cdot K$, 
$Q(g, R)$, $Q(g, S)$, $Q(g, C)$, $Q(g, P)$, $Q(g, W)$, $Q(g, K)$, 
$Q(S, R)$, $Q(S, C)$, $Q(S, P)$, $Q(S, W)$, $Q(S, K)$.\\
\indent 
Using the above presented definitions 
we can prove that 
the following identities 
hold on any semi-Riemannian manifold
$(M,g)$, $n \ge 4$, 
(cf. {\cite[Proposition 1.1] {DGlogH}}): 
$R \cdot K = R \cdot C$ and  
\begin{eqnarray*}
& &
K \cdot S = C \cdot S - \frac{r}{(n-1)(n-2)} Q(g,S) ,\\
& &
K \cdot R = C \cdot R - \frac{r}{(n-1)(n-2)} Q(g,R) ,\\
& &
K \cdot K = C \cdot C - \frac{r}{(n-1)(n-2)} Q(g,C) .
\end{eqnarray*}
Moreover, we also have
$R \cdot W = R \cdot R$, 
$R \cdot C = R \cdot K$, 
$C \cdot R = C \cdot W$, 
$C \cdot C = C \cdot K$, 
$W \cdot R = W \cdot W$, 
$W \cdot C = W \cdot K$,
$K \cdot R = K \cdot W$, 
$K \cdot C = K \cdot K$,
$P \cdot R = P \cdot W$ 
and $P \cdot C = P \cdot K$.
%
\section{\bf Pseudosymmetry type manifolds}

A semi-Riemannian manifold $(M,g)$, $n \geq 3$, is said 
to be an Einstein manifold \cite{Besse} if at every point of $M$ 
its Ricci tensor $S$ is proportional to the metric tensor $g$, 
i.e. $S = \frac{r}{n}\, g$ on $M$.
In particular, if $S$ vanishes on $M$ then it is called Ricci flat. 
We denote by $U_{S}$ the set of all points of $M$ at which 
$S$ is not proportional to $g$, i.e.
$U_S = \left\{ x \in M : S - \frac{r}{n} g \neq 0 \ \mbox{at} \ x \right \}$.\\
\indent As a generalization of Einstein manifold, the notion of quasi-Einstein manifold 
arose during the study of exact solution of Einstein field equations as well as 
during the investigation of quasi-umbilical hypersurfaces of conformally flat spaces. 
For instance, 
FLRW spacetimes are quasi-Einstein spacetimes. 
The semi-Riemannian manifold $(M,g)$, $n \geq 3$, is
said to be a quasi-Einstein manifold 
if 
$\mathrm{rank}\, (S - \alpha\, g) = 1$ 
on $U_S \subset M$, where $\alpha$ is some function on this set
(see, e.g., 
\cite{{ChDGP}, {D0}, {dd}, {D2}, {dd1}, {dd2}, {dd3}, {DaDum01}, {GLOG1}, {SHY09}, {S1}, {S2}, {SRK15}, {S3}}).
In particular, 
if
$\mathrm{rank}\, S = 1$
on $U_S$
then $(M,g)$ is called Ricci-simple \cite{DRVer1}. 
For instance, the G\"{o}del spacetime is a Ricci-simple manifold (see, e.g., \cite{DHJKS}).
The semi-Riemannian manifold $(M,g)$, $n \geq 3$, is
said to be a $2$-quasi-Einstein manifold 
if $\mathrm{rank}\, (S - \alpha\, g) \leq 2$
on $U_S \subset M$ 
and $\mathrm{rank}\, (S - \alpha\, g) = 2$
on some open non-empty subset of $U_S$, 
where $\alpha$ is some function on $U_S$
(see, e.g.,
\cite{{2016_DGHZhyper}, {DGJZ}, {DGP-TV02}, {SKpp}, {SKgrt}, {SKAA17}}).
Such manifolds also are called generalized 
quasi-Einstein manifolds, cf. \cite{{DaDum01}, {DaDum02}} 
and references therein.
It is easy to check that every non-Einstein warped product manifold with 
an $1$-dimensional base and a semi-Riemannian Einsteinian 
$(n-1)$-dimenional fibre is
a quasi-Einstein manifold. Similarly, 
it is easy to check that every non-quasi-Einstein warped product manifold with 
an $2$-dimensional base 
and a semi-Riemannian Einsteinian $(n-2)$-dimenional fibre is
a $2$-quasi-Einstein manifold.
\\ 
\indent
An extension of the class of Einstein semi-Riemannian manifolds $(M, g)$, $n \ge 3$,
also form Ricci-symmetric manifolds \cite{{CA1}, {CA2}, {CA}}, 
i.e. manifolds with parallel Ricci tensor ($\nabla S = 0)$. 
We note that the scalar curvature 
of every Ricci-symmetric manifold is constant. 
In \cite{GRAY} Gray introduced two classes of manifolds lying between 
the class of Ricci-symmetric manifolds and the class of manifolds of constant scalar curvature, viz., 
the class $\mathcal A$ is the class of manifolds which are cyclic Ricci parallel 
$\left((\nabla_X S)(Y,Z)+(\nabla_Y S)(Z,X)+(\nabla_Z S)(X,Y)=0\right)$ 
and the class $\mathcal B$ is the class with Codazi type Ricci tensor 
$\left((\nabla_X S)(Y,Z)=(\nabla_Y S)(X,Z)\right)$. 
Existence of both classes are given in \cite{SB} (see also \cite{DHJKS}). 
Codazzi type Ricci tensor was extensively  studied by verious authors 
\cite{{BER}, {Besse}, {BOUR}, {DER}, {DER1}, {DER2}, {FER}, {SIM}}. 
Another important subclass of the class of Ricci-symmetric manifolds
form locally symmetric manifolds \cite{{CA1}, {CA2}, {CA}}, i.e. manifolds for which we have $\nabla R = 0$. 
This implies the following integrability condition $\mathcal{R}(X,Y) \cdot R = 0$, in short
\be\label{eq1.080808} 
R\cdot R=0 .
\ee
A semi-Riemannian manifold $(M,g)$, $n \ge 3$, is called semisymmetric \cite{CA} if 
(\ref{eq1.080808}) holds on $M$
and a full classification of such manifolds, in the Riemannian case, 
is given by Szab\'{o} \cite{{SZ}, {SZ1}, {SZ2}}.
Further, a semi-Riemannian manifold $(M, g)$, $n \ge 3$, is said to be pseudosymmetric (or, in the sense of Deszcz) 
\cite{{DecuP-TSVer}, {DES}, {DEG}, {DEHAVER}, {HaVer03}, {HV}, {HaVer7}, 
{SDHJK15}, {SKppsnw}, {SKppsn}, {LV1}, {LV2}, {LV3-Foreword}} 
if the tensors $R \cdot R$ and $Q(g,R)$ are linearly dependent at every point of $M$. 
This is equivalent to
\be\label{eq1.3777} 
R\cdot R = L_R Q(g,R)
\ee 
on 
$U_R = \left\{ x \in M : R - \frac{r}{n(n-1)}G\neq 0 \ \mbox{at} \ x \right \}$, where $L_R$ is a function on $U_R$.
Pseudosymmetric manifolds (in the sense of Deszcz \cite{LV1}) 
are also called Deszcz symmetric spaces (see, e.g., \cite{{DecuP-TSVer}, {HaVer7}, {LV_2014}, {LV3-Foreword}}). 
A pseudosymmetric manifold is
called a pseudosymmetric space of constant type if the function 
$L_R$ is constant \cite{{KowSek01}, {KowSek02}}. 
We mention that 
a geometrical interpretation of (\ref{eq1.3777}), in the Rieman\-nian case, 
is given in \cite{HV}.\\ 
\indent We note that pseudosymmetric tensors arose during the study of semisymmetric totally umbilical submanifolds in manifolds 
admitting semisymmetric generalized curvature tensors \cite{{AD1}, {DR1}, {DR7}}. 
For example, every totally umbilical submanifold of a semisymmetric manifold, with parallel mean curvature vector, 
is pseudosymmetric \cite{{AD1}, {ALM}}. The systematic study on pseudosymmetric manifolds was initiated in \cite{AD1}.
We refer to \cite{{DEHAVER}, {HaVer7}} for a wider presentation related to the last statement.
We mention that \cite{DEG} is the first publication, in which  
a semi-Riemannian manifold satisfying (\ref{eq1.3777})
was called the pseudosymmetric manifold.\\ 
\indent The Schwarzschild spacetime, the Kottler spacetime,
the Reissner-Nordstr\"{o}m spacetime 
and the Reissner-Nordstr\"{o}m-de Sitter spacetime
satisfy (\ref{eq1.3777})
with non-zero function $L_R$ \cite{DeVerVra} 
(see also \cite{{P109}, {HaVer03}}). 
We also refer to
\cite{{D0}, {DDVV}, {DeKow}, {DePlSc}, {DeScherf}, {Kow2}}
for further results on pseudosymmetric spacetimes.
For instance, a family of curvature conditions satisfied by
the Reissner-Nordstr\"{o}m-de Sitter spacetime
was determined in \cite{Kow2}. The Schwarzschild spacetime was discovered in 1916 by Schwarzschild
and independently by Droste, during
their study on solutions of Einstein's equations,
see, e.g., 
{\cite[Section 23.3] {Blau}}
and
\cite{Volker} and references therein.
It seems that the Schwarzschild spacetime, 
the Reissner-Nordstr\"{o}m spacetime, as well as some
Friedmann-Lema{\^{\i}}tre-Robertson-Walker (FLRW) spacetimes 
are the ``oldest'' examples of a non-semisymmetric 
pseudosymmetric warped product manifolds (cf. \cite{{DGJZ}, {DEHAVER}}).\\ 
\indent
Pseudosymmetric manifolds form a subclass of the class of Ricci pseudosymmetric manifolds.
A semi-Riemannian manifold $(M,g)$, $n \geq 3$, is said to be Ricci-pseudosymmetric (\cite{{DES}, {DH}}) 
if the tensors $R \cdot S$ and $Q(g,S)$ are linearly dependent at every point of $M$. 
This is equivalent to
\be\label{eq1.3777888} 
R\cdot S = L_{S} Q(g,S)
\ee 
on 
$U_S$, where $L_S$ is a function on this set.
Ricci-pseudosymmetric manifolds are also called Ricci Deszcz symmetric spaces (see, e.g., \cite{LV3-Foreword}).
A Ricci-pseudosymmetric manifold is
called a Ricci-pseudosymmetric manifold of constant type if the function 
$L_S$ is constant 
\cite{Glog2005}. 
It is obvious that
(\ref{eq1.3777}) implies (\ref{eq1.3777888}). The converse statement is not true, provided that $n \geq 4$, 
(see, e.g., \cite{DGHSaw}).
However, (\ref{eq1.3777})  and
(\ref{eq1.3777888}) are equivalent on every $3$-dimensional semi-Riemannian manifold. 
The conditions 
(\ref{eq1.3777}) and (\ref{eq1.3777888})
are equivalent on every $4$-dimensional warped products \cite{D30}.
The conditions  
(\ref{eq1.3777}) and (\ref{eq1.3777888})
are also equivalent on hypersurfaces isometrically immersed 
in $5$-dimensional semi-Riemannian space of constant curvature \cite{1999_DDSVY}.
It is known that
every warped product manifold with 
an $1$-dimensional base and a semi-Riemannian Einsteinian 
$(n-1)$-dimenional fibre is
a Ricci-pseudosymmetric manifold 
\cite{{DGHSaw}, {DH}, {dd3}}.  
For further results on 
Ricci-pseudosymmetric manifolds we refer to \cite{DGHSaw}.
We mention  that 
a geometrical interpretation of
(\ref{eq1.3777888}), in the Riemannian case, 
is given in \cite{JHSV}.\\
\indent
We denote by $U_C$ the set of all points of a semi-Riemannian manifold $(M,g)$, $n \ge 4$, at which 
$C \neq 0$. 
We note that $U_S \cup U_C = U_R$, see, e.g., \cite{DGHHY-2013}.\\
\indent
A semi-Riemannian manifold $(M, g)$, $n \ge 4$, is said to be a manifold with pseudosymmetric Weyl tensor  
\cite{{D30}, {DES}, {DGHHY-2013}, {DEHAVER}, {43}, {46}} if the tensors $C \cdot C$ 
and $Q(g,C)$ are linearly dependent at every point of $M$. 
This is equivalent to
\be\label{eq1.3777999} 
C\cdot C = L_C Q(g,C)
\ee 
on 
$U_C$, where $L_C$ is a function on this set. 
Every warped product manifold with an $2$-dimensional base and 
a $2$-dimenional fibre is a manifold with 
pseudosymmetric Weyl tensor {\cite[Theorem 2] {D30}}. Recently in \cite{DGJZ} it was proved that this 
statement is also true when the fibre is an $(n-2)$-dimensional space of constant curvature, $n \ge 4$.
Thus in particular, the $4$-dimensional spacetime with the metric (\ref{eq1.2}), as well as the $5$-dimensional spacetime with 
the metric (\ref{eq3.4}) are 
$2$-quasi-Einstein manifolds with pseudosymmetric Weyl tensor. 
It may be mentioned that the G\"{o}del spacetime satisfies (\ref{eq1.3777999}) (see, \cite{DHJKS}). 
We refer to
\cite{{DHJKS}, {2019_DHJKS_Erratum}, {SAA17}, {SKpp}, {SK16srs}, {SKAA17}, {SKS17}, {SRK15}}
for examples of various pseudosymmetric type structures.\\
\indent
As it was stated in {\cite[Theorem 3.1]{46}}, if $(M, g)$, $n \ge 4$, 
is a pseudosymmetric manifold with pseudosymmetric Weyl tensor 
then 
\be\label{eq1.37778989} 
Q(S - \alpha  g, C - \beta  G) = 0
\ee 
on 
$U_C$, where $\alpha $ and $\beta $ are some functions on this set. Moreover,
from 
(\ref{eq1.37778989}) it follows that at all points of $U_S \cap U_C$, 
at which 
$\mbox{rank}\, (S - \alpha  g ) > 1$, 
we have (cf. {\cite[Theorem 3.2] {46}}) 
\be\label{eq1.37770909} 
R = L_{1} S \wedge S + L_{2} g \wedge S + L_{3} g \wedge g ,
\ee
where $L_{1}$, $L_{2}$ and $L_{3}$ are some functions on this set.
We refer to  \cite{{DGHSaw}, {DeKow}, {DePlSc}, {DePlSc}, {DeScherf}, {GLOG2}, {Kow1}, {Kow2}}
for results on manifolds satisfying (\ref{eq1.37770909}). 
The manifold satisfying (\ref{eq1.37770909}) is said to be Roter type manifold.
Roter type manifolds are also called Roter spaces.\\
\indent
Some comments on pseudosymmetric manifolds (also called Deszcz symmetric spaces),
as well as Roter spaces, are given in {\cite[Section 1] {DecuP-TSVer}}:
"{\sl{From a geometric point of view, the Deszcz symmetric spaces may well
be considered to be the simplest Riemannian manifolds next to the real space forms.}}" 
and 
"{\sl{From an algebraic point of view, Roter spaces may well be considered to
be the simplest Riemannian manifolds next to the real space forms.}}"
For further comments we refer to \cite{LV3-Foreword}.
Recently Roter spaces admitting geodesic mappings were studied in \cite{DH-2018}.

\section{\bf Geometric structures admitting by interior black hole spacetime}
Let $(B, \overline{g})$ and $(F, \widetilde{g})$ be semi-Riemannian manifolds of dimension $p \geq 1$ and $n-p \geq 1$, respectively, 
covered by the coordinate charts $\{U; x^{a}\}$ and $\{V; y^{\alpha}\}$, respectively. Let $f$ be a smooth positive function on $B$. 
The warped product $M=B\times_{f} F$ is the product manifold $B\times F$ furnished with the metric  
$g=\pi^{*}(g_{B})+(f\circ \pi)\sigma^{*}(g_{F})$, where $\pi$ and $\sigma$ are the projections of $B\times F$ onto $B$ and $F$, 
respectively  
\cite{{Kruchkovich-01}, {Kruchkovich-02}}. 
The manifold $B$ is called the base of $M=B\times F$, and $F$ the fiber. 
We mention that for the warped product manifold $B \times_f F$, the metric can also be considered \cite{{BN}, {Chen-2017}}
as 
$g=\pi^*(g_B)+(f\circ \pi)^2 \sigma^*(g_F)$. However, throughout the paper we will consider the former warped product metric but not later.\\
\indent Let $\{\overline{U}\times \widetilde{V}; x^{1},...,x^{p}, x^{p+1}=y^{1},...,x^{n}=y^{n-p}\}$ 
be a product chart for $B\times F$. The local components of the metric $g=\overline{g}\times_{f} \widetilde{g}$ 
with respect to this chart are given by the following:\\
$g_{hk}=\overline{g}_{ab}$ if $h=a$ and $k=b$, $g_{hk}=f\widetilde{g}_{\alpha \beta}$ 
if $h=\alpha$ and $k=\beta$ and $g_{hk}=0$, otherwise, where $a,b,c,... \in \{1,...,p\}$, 
$\alpha$, $\beta$, $\gamma$, ... $\in$ $\{p+1,...,n\}$ and $h, i, j, k, l, m \in \{1,2,...,n\}$. 
We will mark by bars (resp., by tildes) objects formed from $\overline{g}$ (resp. $\widetilde{g}$).\\
The local components $\Gamma^h_{jk}$ of the Levi-Civita connection $\nabla$ of $B\times_{f} F$ are given by the following:
\be\label{eq:4.2}
\Gamma^a_{bc}=\overline{\Gamma}^a_{bc},\,\,\,\, 
\Gamma^\alpha_{\beta \gamma}=\widetilde{\Gamma}^\alpha_{\beta \gamma},\,\,\,\,\,\,
\Gamma^\alpha_{ab}=\Gamma^a_{\alpha b}=0,
\ee
\be\label{eq:4.3}
\Gamma^a_{\beta \gamma}=-\frac{1}{2}\overline{g}^{ab}f_{b}\widetilde{g}_{\beta \gamma},\,\, 
\Gamma^\alpha_{ a \beta }=\frac{1}{2f}f_{a}\delta^{\alpha}_{\beta}, \,\, 
f_{a}=\partial_{a} f=\frac{\partial f}{\partial x^{a}}.
\ee
The local components
\be\nonumber
R_{hijk}=g_{hl}R^{l}_{ijk}=g_{hl}(\partial_{k}\Gamma^{l}_{ij}-\partial_{j}\Gamma^{l}_{ik}+  \Gamma^{m}_{ij}\Gamma^{l}_{mk}-\Gamma^{m}_{ik}\Gamma^{l}_{mj}),\,\,\, \partial_{k}=\frac{\partial}{\partial x^{k}},
\ee
of the Riemann-Christoffel curvature tensor $R$ and the local components $S_{ij}$ of the Ricci tensor $S$ of the warped product  $\overline{B}\times_{f} \widetilde{F}$ which may not vanish identically are the following:
\be\label{eq:4.4}
R_{abcd}=\overline{R}_{abcd},\,\,\,\, R_{\alpha ab \beta}=-\frac{1}{2} T_{ab}\widetilde{g}_{\alpha \beta},\, \,\,R_{\alpha \beta \gamma \delta}=f\widetilde{R}_{\alpha \beta \gamma \delta}-\frac{\Delta_{1}f}{4}\widetilde{G}_{\alpha \beta \gamma \delta},
\ee
\be\label{eq:4.5}
S_{ab}=\overline{S}_{ab}-\frac{n-p}{2f}T_{ab},\,\,\,\, 
S_{\alpha \beta}=\widetilde{ S}_{\alpha \beta}-\frac{1}{2} \left(tr(T)+\frac{n-p-1}{2f}\Delta_{1}f\right)  \widetilde{g}_{\alpha \beta},
\ee
\be\label{eq:4.6}
T_{ab}=\overline{\nabla}_{b}f_{a}-\frac{1}{2f}f_{a}f_{b},\,\,\,\,tr(T)=\overline{g}^{ab}T_{ab},\,\, \Delta_{1}f= \overline{g}^{ab}f_{a}f_{b},
\ee
and $T$ is the $(0,2)$-tensor with the local components $T_{ab}$. 
The scalar curvature $r$ of $\overline{B}\times_{f} \widetilde{F}$ satisfies the following:
\be\label{eq:4.7}
r=\overline{r}+\frac{\widetilde{r}}{f}-\frac{n-p}{f}\left(tr(T)+\frac{n-p-1}{4f} \Delta_{1}f\right).
\ee
For further details about warped products, we refer to \cite{B}. Warped product pseudosymmetric 
and Ricci pseudosymmetric manifolds are studied 
among others in
\cite{{D0}, {D30}, {DEG}, {DePlSc}, {DeScherf}, {43}}. Also we refer to \cite{k1} 
for the weakly symmetric and weakly Ricci symmetric warped product manifolds.

Let $(M,g)$ be the manifold with the metric $g$ defined by (\ref{eq1.3}). 
Using the above presented formulas we can compute 
the local components of tensors formed by the metric tensor defined by  
(\ref{eq1.3}) (see Section 5). 
We set $\dot{\xi} = \frac{d \xi}{d t}$ and $\ddot{\xi} = \frac{d \dot{\xi }}{d t}$.
We can check that
$S = ( r / 4) g$ 
at all points of $M$ at which $\ddot{\xi} = 2 \dot{\xi} / t$
and 
$\mathrm{rank} ( S +  (2 \dot{\xi} / t^{2} ) g ) = 2$ 
at remains points of $M$, i.e. on $U_{S} \subset M$.
Thus $(M,g)$ is a $2$-quasi-Einstein manifold. 
We can also check that
$S^{2} = (r/2) S - ( 2 \dot{\xi} \ddot{\xi} / t^{3} ) g$,
and in a consequence,
$S^{2} - ( \mathrm{tr}({\mathcal{L}}^{2}) / 4 ) g = (r/2) ( S - ( r / 4) g)$
on $U_{S} \subset M$, where 
$S^{2}$ is a $(0,2)$-tensor with the local components $S^{2}_{ij} = S_{ik}g^{kl}S_{lj}$
and 
$\mathrm{tr}({\mathcal{L}}^{2}) = g^{kl}S^{2}_{kl}$.
Thus the considered manifold is a partially Einstein manifold. 

We recall that 
a semi-Riemannian manifold $(M,g)$, $n \geq 3$, is said to be partially Einstein manifold
{\cite[p. 20] {LV_2014}} if 
$S^{2} = \alpha S + \beta g$ on $U_{S} \subset M$, where $\alpha $ and $\beta$ are some functions on this set.
It is easy to verify that every quasi-Einstein manifold is partially Einstein
(see, e.g., {\cite[Introduction] {GLOG1}}). 
It is obvious that the converse statement is not true. 
We also have \cite{{LV_2014}, {LV3-Foreword}}:
A conformally flat semi-Riemannian manifold of dimension $n \geq 4$
is a Deszcz symmetric space if and only if it is partially Einstein.
 
We present now results related to the interior black hole metric (\ref{eq1.3}). 
\begin{theo}
Interior black hole metric (\ref{eq1.3}) 
satisfies the following:\\
\indent (i) $R \cdot Z = L_1 Q(g, Z)$, \ $L_1 = \frac{\xi - t \dot{\xi}}{t^3}$,\\
\indent (ii) $C \cdot Z = L_2 Q(g, Z)$, \ $L_2 = \frac{6 \xi - 4t \dot{\xi} + t^2 \ddot{\xi}}{6 t^3}$,\\
\indent (iii) $W \cdot Z = L_2 Q(g, Z)$,\\
\indent (iv) $K \cdot Z = L_3 Q(g, Z)$, \ $L_3 = \frac{2\xi + t^2 \ddot{\xi}}{2t^3}$,\\
\indent (v) $P \cdot S = L_1 Q(g, S)$,\\
\indent (vi) $P \cdot Z = L_1 Q(g, Z) - \frac{1}{3}Q(S, Z)$,\\
where $Z$ is any one of $R, S, C, W, K$ and $P$.
\end{theo}
From above theorem it follows that $(a)$ if $\xi = C_1 t$ then $R \cdot Z = 0$ and $P \cdot S = 0$, $(b)$ if $\xi = C_1 t^2 + C_2 t^3$ 
then $C \cdot Z = 0$ and $W \cdot Z = 0$, and $(c)$ if 
$\xi = \sqrt{t} \left( C_2 \cos\left(\frac{\sqrt{7}}{2}\log t\right) + C_1 \sin\left(\frac{\sqrt{7}}{2}\log t \right) \right)$ 
then $K \cdot Z = 0$, where $C_1$ and $C_2$ are some constants. 
Further, we have
\begin{theo}
Interior black hole metric (\ref{eq1.3}) satisfies the following:\\
\indent (i) $C \cdot K = W \cdot K$, \ \ \ \ \ $C \cdot K = W \cdot C$, \ \ \ \ \ $C \cdot C = W \cdot C$,\\
\indent (ii) $W \cdot K = W \cdot C$, \ \ \ \ \ $C \cdot K = C \cdot C$ \ ((ii) follows from (i)),\\
\indent (iii) $W \cdot K = C \cdot C$ \ ((iii) follows from (i) and (ii)),\\
\indent (iv) $C \cdot W = W \cdot R$, \\
\indent (v) $R \cdot S = P \cdot S$, \ \ \ \ \ $C \cdot S = W \cdot S$,\\
\indent (vi) $L_3 \ R \cdot K = L_1 \ K \cdot C$,\\
\indent (vii) $R \cdot W - W \cdot R = L_5 \ Q(g, R)$, \ $L_5 = - \frac{2\dot{\xi} + t \ddot{\xi}}{6t^2} = \frac{r}{12}$,\\
\indent (viii) $C \cdot K - K \cdot C = L_6 \ Q(g, C)$, \ $L_6 = - \frac{2 \dot{\xi} + t \ddot{\xi}}{3t^2}  = \frac{r}{6}$,\\
\indent (ix) $C \cdot R - Q(S, C) = L_7 \ Q(g, C)$, \ $L_7 = \frac{3 \xi + t \dot{\xi} + 2t^2 \ddot{\xi}}{3t^3}$,\\
\indent (x) $R \cdot R - Q(S, R) = L_8 \ Q(g, C)$, \ $L_8 = \frac{6\xi^2 - 4 t \xi \dot{\xi} - 2t^2 \dot{\xi}^2 + 4t^2 \xi \ddot{\xi}}{t^3(6 \xi - 4 t \dot{\xi} + t^2 \ddot{\xi})}$,\\
\indent (xi) $L_6 L_9 \ R \cdot W + L_1 L_{10} \ Q(S, W) = L_1 L_{11} \ Q(S, R)$, \ $L_9 = - 3t (6\dot{\xi}(\xi- t\dot{\xi})+t(t\dot{\xi}+3\xi)\ddot{\xi})$,\\
\indent \ \ \ \ \ \ $L_{10} = 6(-t^2 \dot{\xi}^2+2 t \xi(t \ddot{\xi}-\dot{\xi})+3 \xi^2)$, \ $L_{11} = t^2(-t^2 \ddot{\xi}^2+2 \dot{\xi}^2+2t\dot{\xi} \ddot{\xi})+6t\xi(t \ddot{\xi}-4 \dot{\xi})+18 \xi^2$,\\
\indent (xii) $L_3 \ R \cdot K + L_{12} \ K \cdot R + L_3 \ Q(S, K)=0$, \ $L_{12} = \frac{-2\xi + 4 t \dot{\xi} + t^2 \ddot{\xi}}{2t^3}$,\\
\indent (xiii) $L_3 L_{13} \ R \cdot K - L^2_1 L_2 \ K \cdot R + L_1 L_2 L_3 \ Q(S, R) = 0$, \ $L_{13} = \frac{3\xi^2 - 2 t \xi \dot{\xi} - t^2 \dot{\xi}^2 + 2t^2 \xi \ddot{\xi}}{3 t^6}$,\\
\indent (xiv) $L_2 \ C \cdot W - L_7 \ W \cdot C = L_2 \ Q(S, C)$,\\
\indent (xv) $L_2 L_{14} \ C \cdot W - \frac{1}{18t^6}L_{11} \ W \cdot C = L^2_2 \ Q(S, W)$, \ \ $L_{14} = \frac{3\xi - 5t \dot{\xi} - t^2 \ddot{\xi}}{3t^3}$,\\
\indent (xvi) $L_{15} \ C \cdot K + L^2_2 \ Q(S, K) = L_2 L_{16} \ Q(S, C)$, \ $L_{15} = \frac{4t^2\dot{\xi}+4t^3\dot{\xi}\ddot{\xi}+t^4 \ddot{\xi}^2)}{3t^6}$, \ \ \ $L_{16} = \frac{2\xi -4t \dot{\xi}-t^2 \ddot{\xi}}{2t^3}$,\\
\indent (xvii) $18t^6 L_3 L_{11} \ W \cdot K - L^2_2 L_{14} \ K \cdot W + L^2_2 L_3 \ Q(S, W) = 0$,\\
\indent (xviii) $L_1 L_3 \ W \cdot K + L_2 L_{12} \ K \cdot W + L_2 L_3 \ Q(S, K)=0$,\\
\indent (xix) Roter type condition with $R = \frac{\phi}{2}S\wedge S + \mu g\wedge S + \eta G$, where 
$$
\phi = - \frac{ 6 t \xi - 4 t^2 \dot{\xi} + t^3 \ddot{\xi}}{(t \ddot{\xi} - 2 \dot{\xi})^2}, 
\mu = -\frac{6 \xi \dot{\xi} - 6 t \dot{\xi}^2 + 3 t \xi \ddot{\xi} + t^2 \xi \ddot{\xi}}{t(t \ddot{\xi} - 2 \dot{\xi})^2}, 
\eta = -\frac{2(4 \xi \dot{\xi}^2 - 4 t \dot{\xi}^3 + 2 t \xi \dot{\xi} \ddot{\xi} + t^2 \xi \ddot{\xi}^2)}{t^3(t \ddot{\xi} - 2 \dot{\xi})^2}.
$$
\end{theo}
From above theorem it follows that 
$(a)$ if $\xi = C_1 t$ 
then $R \cdot K = K \cdot C$, $(b)$ if $\xi = -\frac{C_1}{t} + C_2$ 
then $R \cdot W = W \cdot R$ and $C \cdot K = K \cdot C$, $(c)$ 
if $\xi = t^{\frac{1}{4}} \left( C_2 \cos\left(\frac{\sqrt{23}}{4}\log t\right) 
+ C_1 \sin\left(\frac{\sqrt{23}}{4}\log t\right) \right)$ 
then $C \cdot R = Q(S, C)$,
where $C_1$ and $C_2$ are some constants.\\
\indent In \cite{DePlSc} it was proved that any Roter type manifold satisfies the relation (x) with 
$L_8 = L_1 + \phi^{-1}\mu = (n-2)\phi^{-1}(\mu^2 - \phi\eta)$, $n = 4$, and the converse is also true as follows from (x) and (xix).\\

\indent 
Let $(M,g)$ be the manifold with the interior black hole metric in $5$-dimension is given by \cite{{HK1}, {HK}, {HK2}}
\be\label{eq3.4}
ds^2 = - \left(\frac{2\xi}{t^2} - 1\right)^{-1}dt^2 + \left(\frac{2\xi}{t^2} - 1\right)dz^2 + t^2\left(d\theta^2 + \sin^2\theta d\phi^2 + \sin^2 \theta \sin^2 \phi d\psi^2\right),
\ee
where $\xi$ is the function of time. 
Similarly as above, we can obtain curvature properties of the interior 
black hole metric (\ref{eq3.4}). In particular, we again set
$\dot{\xi} = \frac{d \xi}{d t}$ and $\ddot{\xi} = \frac{d \dot{\xi }}{d t}$.
We can check that
$S = ( r / 5) g$ 
at all points of $M$ at which $\ddot{\xi} = 3 \dot{\xi} / t$ 
and 
$\mathrm{rank} ( S +  (2 \dot{\xi } / t^{3}) g ) = 2$
at remains points of $M$, i.e. on $U_{S} \subset M$.
Thus $(M,g)$ is a $2$-quasi-Einstein manifold. 
We can also check that
$S^{2} - ( \mathrm{tr}({\mathcal{L}}^{2}) / 5 ) g = ( (r/2) + (\dot{\xi } / t^{3})  ) ( S - ( r / 5) g)$
on $U_{S} \subset M$.
Thus the considered manifold is a partially Einstein manifold. 
 
Further, we have the following curvature properties of the considerd metric (\ref{eq3.4}). 
\begin{theo}
$5$-dimensional interior black hole metric 
(\ref{eq3.4})
satisfies the following:\\
\indent (i) $R \cdot Z = N_1 Q(g, Z)$, \ $N_1 = \frac{2\xi - t \dot{\xi}}{t^4}$,\\
\indent (ii) $C \cdot Z = N_2 Q(g, Z)$, \ $N_2 = \frac{12 \xi - 6t \dot{\xi} + t^2 \ddot{\xi}}{6 t^4}$,\\
\indent (iii) $W \cdot Z = N_3 Q(g, Z)$, \ $N_3 = \frac{20 \xi - 8t \dot{\xi} + t^2 \ddot{\xi}}{10 t^4}$,\\
\indent (iv) $K \cdot Z = N_4 Q(g, Z)$, \ $N_4 = \frac{6\xi - 2 t \dot{\xi} + t^2 \ddot{\xi}}{3t^4}$,\\
\indent (v) $P \cdot S = N_1 Q(g, S)$,\\
\indent (vi) $P \cdot Z = N_1 Q(g, Z) - \frac{1}{4}Q(S, Z)$,\\
where $Z$ is any one of $R, S, C, W, K$ and $P$.
\end{theo}
From above theorem it follows that $(a)$ if $\xi = C_1 t^2$ then $R \cdot Z = 0$ and $P \cdot S = 0$, $(b)$ if $\xi = C_1 t^3 + C_2 t^4$ 
then $C \cdot Z = 0$, $(c)$ if $\xi = C_1 t^5 + C_2 t^4$ then $W \cdot Z = 0$, 
and $(d)$ if 
$\xi = t^{\frac{3}{2}} \left( C_2 \cos\left(\frac{\sqrt{15}}{2}\log t\right) 
+ C_1 \sin\left(\frac{\sqrt{15}}{2}\log t\right) \right)$ 
then $K \cdot Z = 0$, where $C_1$ and $C_2$ are some constants.
\begin{theo}
$5$-dimensional interior black hole metric (\ref{eq3.4}) satisfies the following:\\
\indent (i) $R \cdot S = P \cdot S$,\\
\indent (ii) $R \cdot W - W \cdot R = N_5 Q(g, R)$, \ $N_5 = - \frac{2\dot{\xi} + t \ddot{\xi}}{10t^3}$,\\
\indent (iii) $C \cdot K - K \cdot C = N_6 Q(g, C)$, \ $N_6 = - \frac{2 \dot{\xi} + t \ddot{\xi}}{6t^3}$,\\
\indent (iv) $C \cdot R - Q(S, C) = N_7 Q(g, C)$, \ $N_7 = \frac{4 \xi + t^2 \ddot{\xi}}{2t^4}$,\\
\indent (v) $R \cdot R - Q(S, R)= N_8 \ Q(g, C)$, \ $N_8 = \frac{3(8\xi^2 - 4 t \xi \dot{\xi} - t^2 \dot{\xi}^2 + 2t^2 \xi \ddot{\xi})}{t^4(12\xi - 6t\dot{\xi} + t^2 \ddot{\xi})}$,\\
\indent (vi) $N_5 N_9 \ R \cdot W + N_1 N_{10} \ Q(S, R) + N_1 N_8 \ Q(S, W) = 0$, \ \ $N_9 = -\frac{3\dot{\xi}(4\xi - 3t \dot{\xi})+ t\ddot{\xi}(4\xi+t\dot{\xi})}{2t^7}$,\\
\indent \ \ \ \ \ \ \ \ \ \ \ \ \ $N_{10}=\frac{-80\xi^2 + 64t\xi \dot{\xi} - 2t^2\dot{\xi}^2 -4t^2 \ddot{\xi}(2\xi + t \dot{\xi})+ t^4 \ddot{\xi}^2}{20t^8}$,\\
\indent (vii) $N_4 \ R \cdot K + N_{11} \ K \cdot R + N_4 \ Q(S, K) = 0$, \ $N_{11} = \frac{-6\xi + 6t \dot{\xi} + t^2 \ddot{\xi}}{3t^4}$,\\
\indent (viii) $N_4 N_8 \ R \cdot K - N^2_1 N_2 \ K \cdot R + N_1 N_2 N_4 \ Q(S, R) = 0$,\\
\indent (ix) $N_3 \ C \cdot W - N_7 \ W \cdot C = N_3 \ Q(S, C)$,\\
\indent (x) $N_3 N_{12} \ C \cdot W + N_{10} \ W \cdot C = N_3 N_2 \ Q(S, W)$, \ \ $N_{12}=\frac{20 \xi - 16 t \dot{\xi} - 3t^2 \ddot{\xi}}{10t^4}$,\\
\indent (xi) $N_{13} \ C \cdot K + N_2 N_{11} \ Q(S, C) + N^2_2 \ Q(S, K) = 0$, \ \ $N_{13} = \frac{6\dot{\xi}^2 + 5t \dot{\xi}\ddot{\xi}+t^2 \ddot{\xi}^2}{6t^6}$,\\
\indent (xii) $N_4 N_{10} \ W \cdot K + N_{12} N_{14} \ K \cdot W = N_4 N_{14} \ Q(S, W)$,\\
\indent \ \ \ \ \ \ \ \ \ \ \ \ \ \ $N_{14}=\frac{240\xi^2 - 216t \xi \dot{\xi} + 48t^2 \dot{\xi}^2 + 32t^2 \xi \ddot{\xi} - 14t^3 \dot{\xi} \ddot{\xi} + t^4 \ddot{\xi}^2}{60t^8}$,\\
\indent (xiii) $N_1 N_4 \ W \cdot K + N_3 N_{11} \ K \cdot W + N_3 N_4 \ Q(S, K) = 0$,\\
\indent (xiv) Roter type condition with $R = \frac{\phi}{2}S\wedge S + \mu g\wedge S + \eta G$, where 
\begin{eqnarray*}
& &
\phi  =  -\frac{t^2( 12\xi + t(t \ddot{\xi} - 6 \dot{\xi}))}{(3 \dot{\xi} - t \ddot{\xi})^2}, \ \  
\mu   =  -\frac{3 \dot{\xi}(4\xi - 3 t \dot{\xi}) + t (4 \xi + t \dot{\xi})\ddot{\xi}}{t(3 \dot{\xi} - t \ddot{\xi})^2},\\ 
& &
\eta  =  -\frac{2(-6t\dot{\xi}^3 + \xi( 9\dot{\xi}^2 + 2t \dot{\xi} \ddot{\xi} + t^2 \ddot{\xi}^2))}{t^4(3 \dot{\xi} - t \ddot{\xi})^2}.
\end{eqnarray*}
\end{theo}
From above theorem it follows that $(a)$ if $\xi = -\frac{C_1}{t} + C_2$ then $R \cdot W = W \cdot R$ and $C \cdot K = K \cdot C$, 
$(b)$ 
if $\xi = \sqrt{t} \left( C_2 \cos\left(\frac{\sqrt{15}}{2}\log t\right) \
+ C_1 \sin\left(\frac{\sqrt{15}}{2}\log t\right) \right)$ then $C \cdot R = Q(S, C)$,
and $(c)$ if $\xi = \frac{C_2 \left(3 + 10 t^5 C_1\right)^{\frac{2}{3}}}{t^2}$ then $R \cdot R = Q(S, R)$,
where $C_1$ and $C_2$ are some constants.\\
However, we note that for the $5$-dimensional interior black hole metric
the following tensors are non-zero tensors\\
\indent (i) $C \cdot K - W \cdot K$, (ii) $C \cdot K - W \cdot C$, (iii) $C \cdot C - W \cdot C$, (iv) $C \cdot W - W \cdot R$, 
(v) $C \cdot S - W \cdot S$ and (vi) $R \cdot K - K \cdot C$.\\
\indent In \cite{DePlSc} it was proved that any Roter type manifold 
satisfies the relation (v) with 
$N_8 = N_1 + \phi^{-1}\mu = (n-2)\phi^{-1}(\mu^2 - \phi\eta)$, $n=5$, 
and hence from (xiv), it follows that the converse of the result is also true.\\

We can check that 
the $4$-dimensional spacetime with the metric (\ref{eq1.3}) 
and the $5$-dimensional spacetime with the metric (\ref{eq3.4}) 
are non-quasi Einstein manifolds. Furthermore
from the considerations presented in Section 3 and 
Theorems 4.1 (i) and 4.4 (ii) it follow that
those spacetimes
satisfy (\ref{eq1.37770909}) 
and hence are Roter type spacetimes. \\
\indent We also note that interior black hole metrics 
(\ref{eq1.3}) and (\ref{eq3.4})
do not admit any one of the following structures:
Ricci semisymmetric, 
quasi-Einstein, 
Codazzi type Ricci tensor, 
cyclic Ricci symmetric, 
Chaki pseudo symmetric, 
Chaki pseudo Ricci symmetric, 
weakly symmetric, 
weakly Ricci symmetric, 
hyper generalized recurrent, 
weakly generalized recurrent, 
quasi generalized recurrent, 
as well as 
any pseudosymmetric type structure defined by $R \cdot Z = L Q(S, Z)$,  $C \cdot Z = L Q(S, Z)$, $W \cdot Z = L Q(S, Z)$, 
$K \cdot Z = L Q(S, Z)$, $P \cdot Z = L Q(S, Z)$, and $P \cdot Z = L Q(g, Z)$, respectively,
where $L$ is any smooth function and $Z$ is any one of the tensors $R, C, W, K$ and $P$. 
It can also be mentioned that both the interior black hole metrics does not realize 
any one of the generalized Einstein metric condition (i) $R\cdot C - C\cdot R = L_1 Q(g,R)$, 
(ii) $R\cdot C - C\cdot R = L_2 Q(g,C)$, (iii) $R\cdot C - C\cdot R = L_3 Q(S,R)$ 
and (iv) $R\cdot C - C\cdot R = L_4 Q(S,C)$. 
We mention that a survey on semi-Riemannian manifolds satisfying the last four conditions is given in \cite{DGHSaw}.
\section{\bf Appendix}

\noindent
{\bf{Part I.}}
From 
(\ref{eq1.3}) 
and (\ref{eq:4.2})-(\ref{eq:4.7}), 
the local components of the Christoffel symbols of second kind,
the curvature tensor and the Ricci tensor (upto symmetry)
which may not vanish identically are the following:
$$\Gamma^{1}_{11} = -\frac{\xi-t \dot{\xi}}{t^2-2 t\xi} = -\Gamma^{2}_{12}, \ \ 
\Gamma^{3}_{13} = \Gamma^{4}_{14} = \frac{1}{t}, \ \ \Gamma^{1}_{22} = -\frac{(t-2 \xi) (t \dot{\xi}-\xi)}{t^3},$$
$$\Gamma^{1}_{33} = 2 \xi-t, \ \ \Gamma^{4}_{34} = \cot \theta, \ \ \Gamma^{1}_{44} = -(t-2 \xi) \sin^2\theta, \ \ 
\Gamma^{3}_{44} = -\sin \theta \cos \theta,$$
$$ R_{1212} = -\frac{t^2 \ddot{\xi}-2 t \dot{\xi}+2 \xi}{t^3}, \ \ \ R_{1313} = \frac{t \dot{\xi}-\xi}{t-2 \xi}, \ \ \ 
R_{1414} = \frac{\sin^2\theta (t\dot{\xi}-\xi)}{t-2 \xi},$$
$$ R_{2323} = -\frac{(t-2 \xi) (t \dot{\xi}-\xi)}{t^2}, \ \ \ R_{2424} = -\frac{(t-2 \xi) \sin^2\theta (t \dot{\xi}-\xi)}{t^2}, \ \ \ 
R_{3434} = 2 t \xi \sin^2\theta,$$
\begin{eqnarray*} 
& &
S_{11} = -\frac{\ddot{\xi}}{t-2 \xi}       = - \frac{ \ddot{\xi} }{t} g_{11},\ \ \
S_{22} = \frac{(t-2 \xi) \ddot{\xi}}{t^2}  = - \frac{ \ddot{\xi} }{t} g_{22},\\
& &
S_{33} = -2 \dot{\xi} =  - \frac{2 \dot{\xi}}{t^{2}} g_{33},\ \ \ 
S_{44} = -2 \dot{\xi} \sin^{2} \theta  =  - \frac{2 \dot{\xi}}{t^{2}} g_{44},\ \ \
r = -\frac{2(2\dot{\xi} + t\ddot{\xi})}{t^2} , 
\end{eqnarray*}
where $\dot{\xi}$ denotes the differentiation of $\xi$ with respect to $t$
and 
$\ddot{\xi}$  the second order differentiation of $\xi$ with respect to $t$.
The local components of the Weyl conformal curvature tensor (upto symmetry) which may not vanish identically are given by
$$ C_{1212} = -\frac{6\xi -t(4 \dot{\xi} - t \ddot{\xi})}{3 t^3}, \ \ \ 
C_{1313} = -\frac{6\xi -t(4 \dot{\xi} - t \ddot{\xi})}{6(t-2 \xi)}, \ \ \ 
C_{1414} = -\frac{\sin^2\theta [6\xi -t(4 \dot{\xi} - t \ddot{\xi})]}{6(t - 2 \xi)},$$
$$ C_{2323} = \frac{(t-2\xi) [6\xi -t(4 \dot{\xi} - t \ddot{\xi})]}{6(t-2 \xi)}, \ \ \ 
C_{2424} = \frac{(t-2\xi) \sin^2\theta [6\xi -t(4 \dot{\xi} - t \ddot{\xi})]}{6(t-2 \xi)},$$
$$C_{3434} = -\frac{1}{3} t \sin^2\theta [6\xi -t(4 \dot{\xi} - t \ddot{\xi})].$$
The local components of the covariant derivatives of curvature tensor and Ricci tensor (upto symmetry) 
which may not vanish identically 
are given by:
$$ R_{1212,1}=\frac{- t^3 \xi^3 +3 t^2 \ddot{\xi}-6 t \dot{\xi}+6 \xi}{t^4}, \ \ \ R_{1223,3}=\frac{(t-2 \xi) \left(t (t \ddot{\xi}-3 \dot{\xi})+3 \xi\right)}{t^3} = - R_{2323,1},$$
$$ R_{1224,4}=\frac{(t-2 \xi) \sin^2\theta \left(t (t \ddot{\xi}-3 \dot{\xi})+3 \xi\right)}{t^3} = - R_{2424,1}, \ \ \ R_{1313,1}=\frac{t^2 \ddot{\xi}-3 t \dot{\xi}+3 \xi}{t^2-2 t\xi},$$
$$ R_{1334,4}=\sin^2\theta (3 \xi-t \dot{\xi}) = - R_{1434,3} = -\frac{1}{2}R_{3434,1}, \ \ \ R_{1414,1}=\frac{\sin^2\theta \left(t (t \ddot{\xi}-3 \dot{\xi})+3 \xi\right)}{t (t-2 \xi)},$$
$$ S_{11,1}=\frac{\ddot{\xi}-t \xi^3}{t^2-2 t\xi}, \ \ S_{13,3}=\frac{2 \dot{\xi}}{t}-\ddot{\xi} = \frac{1}{2}S_{33,1},$$
$$ S_{22,1}=\frac{(t-2 \xi) (t \xi^3-\ddot{\xi})}{t^3}, \ \ S_{14,4}=\frac{\sin^2\theta (2 \dot{\xi}-t \ddot{\xi})}{t} = \frac{1}{2}S_{44,1}.$$

In terms of local coordinate system, the local components
$Q(A,T)_{i_{_1} i_{_2} \ldots i_{_k} u v}$ of 
the Tachibana tensor $Q(A,T)$ of an $(0,2)$-tensor $A$ and an $(0,k)$-tensor $T$ are given by
\begin{eqnarray*}
Q(A,T)_{i_{_1} i_{_2} \ldots i_{_k} u v} &=& A_{i_{_1} u} T_{v i_{_2} \ldots i_{_k}} + A_{i_{_2} u} T_{i_{_1} v \cdots i_{_k}} 
+ \cdots + A_{i_{_k} u} T_{i_{_1} i_{_2} \cdots v}\\
&& -A_{i_{_1} v} T_{u i_{_2} \cdots i_{_k}} - A_{i_{_2} v} T_{i_{_1} u \cdots i_{_k}} 
- \cdots - A_{i_{_k} v} T_{i_{_1} i_{_2} \cdots u}.
\end{eqnarray*}
In particular, for a symmetric $(0,2)$-tensor $A$ and a generalized curvature tensor $T$, we have
\begin{eqnarray*}
& & Q(A,T)_{hijklm} = Q(A,T)_{jkhilm} = -  Q(A,T)_{ihjklm} = - Q(A,T)_{hijkml}, \\
& & Q(A,T)_{hijklm} = Q(A,T)_{ijhklm} = Q(A,T)_{jhiklm} , \\
& & Q(A,T)_{hijklm} + Q(A,T)_{jklmhi} + Q(A,T)_{lmhijk} = 0.   
\end{eqnarray*}
If $A$ and $B$ are symmetric $(0,2)$-tensors then 
\begin{eqnarray*}
& & Q(A,B)_{hijk} = Q(A,B)_{ihjk} = - Q(A,B)_{hikj} .
\end{eqnarray*}

If $\mathcal{D}(X,Y) = \mathcal{R}(X,Y)$ then (\ref{ror}) yields
\begin{eqnarray*}
(R\cdot T)_{i_{_1} i_{_2} \cdots i_{_k} u v} 
&=& -g^{p q}(T_{p i_{_2} \cdots i_{_k}}R_{u v q i_{_1}}  +  T_{i_{_1} p \cdots i_{_k}} R_{u v q i_{_2}}
+ \cdots + T_{i_{_1} i_{_2} \cdots p} R_{u v q i_{_k}} ),
\end{eqnarray*}
where $g^{p q}$, $R_{hijk}$ and $T_{i_{_1} i_{_2} \cdots i_{_k}}$ are the local components of the tensors
$g^{-1}$, $R$ and $T$, respectively. 
Similarly in terms of local coordinate system we can write the components of 
$C\cdot T$, $P\cdot T$, $W\cdot T$ and $K\cdot T$. 
Moreover, if $B_{hijk}$ and $T_{hijk}$ are the local components 
of generalized curvature tensors $B$ and $T$ then the local components of the $(0,6)$-tensor $B \cdot T$
are following
\begin{eqnarray*}
(B \cdot T)_{hijklm}  
&=& -g^{p q}(T_{pijk} B_{lmhq}  + T_{hpjk} B_{lmiq}
+ T_{hipk} B_{lmjq} + T_{hijp} B_{lmkq}  ) \\
&=& g^{p q}(T_{p ijk} B_{lmqh}  -  T_{phjk} B_{lmqi}
+ T_{pk hi} B_{lm q j} -  T_{pj hi} B_{lm q k}  ) .
\end{eqnarray*}
We have
\begin{eqnarray*}
& & (B \cdot T))_{hijklm} = (B \cdot T)_{jkhilm} = -  (B \cdot T)_{ihjklm} = - (B \cdot T)_{hijkml}, \\
& & (B \cdot T)_{hijklm} +(B \cdot T)_{ijhklm} + (B \cdot T)_{jhiklm} = 0 .
\end{eqnarray*}

For the tensors $R \cdot R$, $R \cdot S$, $R \cdot C$ and $R \cdot P$ 
we have the following relations:
$$(R \cdot R)_{122313}=\frac{(t \dot{\xi}-\xi) \left(t (t \ddot{\xi}-3 \dot{\xi})+3 \xi\right)}{t^4} = -(R \cdot R)_{121323},$$
$$(R \cdot R)_{143413}=\frac{\sin^2\theta (t \dot{\xi}-3 \xi)(t \dot{\xi}-\xi)}{t (t-2\xi)} = -(R \cdot R)_{133414},$$
$$(R \cdot R)_{122414}=\frac{\sin^2\theta (t \dot{\xi}-\xi)\left(t(t \ddot{\xi}-3\dot{\xi})+3 \xi\right)}{t^4} = (R \cdot R)_{121424},$$
$$-(R \cdot R)_{243423}=\frac{(t-2 \xi) \sin^2\theta (t \dot{\xi}-3 \xi) (t \dot{\xi}-\xi)}{t^3} = (R \cdot R)_{233424};$$
$$(R \cdot S)_{1313}=\frac{(\xi-t \dot{\xi}) (t \ddot{\xi}-2 \dot{\xi})}{t^2 (t-2 \xi)}, \ \ (R \cdot S)_{1414}=\frac{\sin^2\theta (\xi-t \dot{\xi}) (t \ddot{\xi}-2 \dot{\xi})}{t^2 (t-2 \xi)},$$
$$(R \cdot S)_{2323}=\frac{(t-2 \xi) (t \dot{\xi}-\xi) (t \ddot{\xi}-2 \dot{\xi})}{t^4}, \ \ (R \cdot S)_{2424}=\frac{(t-2 \xi) \sin^2\theta (t \dot{\xi}-\xi) (t \ddot{\xi}-2 \dot{\xi})}{t^4};$$
$$(R \cdot C)_{122313}=\frac{(t \dot{\xi}-\xi) \left(t (t \ddot{\xi}-4 \dot{\xi})+6 \xi\right)}{2 t^4} = -(R \cdot C)_{121323},$$
$$(R \cdot C)_{143413}
=\frac{\sin^2\theta (\xi-t \dot{\xi}) \left(t (t \ddot{\xi}-4 \dot{\xi})+6 \xi\right)}{2 t (t-2 \xi)}=-(R \cdot C)_{133414},$$
$$(R \cdot C)_{122414}=\frac{\sin^2\theta (t \dot{\xi}-\xi) \left(t (t \ddot{\xi}-4 \dot{\xi})+6 \xi\right)}{2 t^4} 
=- (R \cdot C)_{121424},$$
$$(R \cdot C)_{243423}=\frac{(t-2 \xi) \sin^2\theta (t \dot{\xi}-\xi) \left(t (t \ddot{\xi}-4 \dot{\xi})+6 \xi\right)}{2 t^3}
=-(R \cdot C)_{233424};$$

$$(R \cdot P)_{122313} = - (R \cdot P)_{131223} = - (R \cdot P)_{232113} = - (R \cdot P)_{121323}=\frac{(t \dot{\xi}-\xi) \left(t (t \ddot{\xi}-3 \dot{\xi})+3 \xi\right)}{t^4},$$
$$(R \cdot P)_{123213} = - (R \cdot P)_{132123} = - (R \cdot P)_{231213} = - (R \cdot P)_{123123}=\frac{(\xi-t \dot{\xi}) \left(t (2 t \ddot{\xi}-7 \dot{\xi})+9 \xi\right)}{3 t^4},$$
$$(R \cdot P)_{131113}=\frac{(t \dot{\xi}-\xi) (t \ddot{\xi}-2 \dot{\xi})}{3 t (t-2 \xi)^2}, \ \ (R \cdot P)_{133313}=\frac{(\xi-t \dot{\xi}) (t \ddot{\xi}-2 \dot{\xi})}{3 (t-2 \xi)},$$
$$(R \cdot P)_{143413} = - (R \cdot P)_{341314} = (R \cdot P)_{134314} = - (R \cdot P)_{341413}=\frac{\sin^2 \theta (\xi-t \dot{\xi}) \left(t (t \ddot{\xi}-5 \dot{\xi})+9 \xi\right)}{3 t (t-2 \xi)},$$
$$(R \cdot P)_{144313} = (R \cdot P)_{344113} = (R \cdot P)_{133414} = - (R \cdot P)_{343114}=-\frac{\sin^2 \theta (t \dot{\xi}-3 \xi) (t \dot{\xi}-\xi)}{t (t-2 \xi)},$$
$$(R \cdot P)_{122414} = - (R \cdot P)_{242114} = - (R \cdot P)_{121424} = - (R \cdot P)_{141224}=\frac{\sin^2 \theta (t \dot{\xi}-\xi) \left(t (t \ddot{\xi}-3 \dot{\xi})+3 \xi\right)}{t^4},$$
$$(R \cdot P)_{124214} = - (R \cdot P)_{241214} = - (R \cdot P)_{124214} = - (R \cdot P)_{142124}=\frac{\sin^2 \theta (\xi-t \dot{\xi}) \left(t (2 t \ddot{\xi}-7 \dot{\xi})+9 \xi\right)}{3 t^4},$$
$$(R \cdot P)_{141114}=\frac{\sin^2 \theta (t \dot{\xi}-\xi) (t \ddot{\xi}-2 \dot{\xi})}{3 t (t-2 \xi)^2}, \ \ (R \cdot P)_{144414}=\frac{\sin^4 \theta (\xi-t \dot{\xi}) (t \ddot{\xi}-2 \dot{\xi})}{3 (t-2 \xi)},$$
$$(R \cdot P)_{232223}=\frac{(t-2 \xi)^2 (t \dot{\xi}-\xi) (t \ddot{\xi}-2 \dot{\xi})}{3 t^5}, \ \ (R \cdot P)_{233323}=\frac{(t-2 \xi) (t \dot{\xi}-\xi) (t \ddot{\xi}-2 \dot{\xi})}{3 t^2},$$
$$(R \cdot P)_{243423} = (R \cdot P)_{342423} = (R \cdot P)_{234324}$$
$$ = - (R \cdot P)_{342324}=\frac{(t-2 \xi) \sin^2 \theta (t \dot{\xi}-\xi) \left(t (t \ddot{\xi}-5 \dot{\xi})+9 \xi\right)}{3 t^3},$$
$$(R \cdot P)_{244323} = (R \cdot P)_{344223} = (R \cdot P)_{233424} = - (R \cdot P)_{343224}=\frac{(t-2 \xi) \sin^2 \theta (t \dot{\xi}-3 \xi) (t \dot{\xi}-\xi)}{t^3},$$
$$(R \cdot P)_{242224}=\frac{(t-2 \xi)^2 \sin^2 \theta (t \dot{\xi}-\xi) (t \ddot{\xi}-2 \dot{\xi})}{3 t^5}, \ \ (R \cdot P)_{244424}=\frac{(t-2 \xi) \sin^4 \theta (t \dot{\xi}-\xi) (t \ddot{\xi}-2 \dot{\xi})}{3 t^2}.$$

For the tensors $Q(g, R)$, $Q(g, S)$, $Q(g, C)$ and $Q(g, P)$ we have the following relations:
$$Q(g,R)_{122313}=-t \ddot{\xi}+3 \dot{\xi}-\frac{3 \xi}{t} = -Q(g,R)_{121323},$$
$$Q(g,R)_{143413}=\frac{t^2 \sin^2\theta (3 \xi-t \dot{\xi})}{t-2 \xi} = -Q(g,R)_{133414},$$
$$Q(g,R)_{121424}=\frac{\sin^2\theta \left(t (t \ddot{\xi}-3 \dot{\xi})+3 \xi\right)}{t}= - Q(g,R)_{122414},$$
$$Q(g,R)_{243423}=(t-2 \xi) \sin^2\theta (t \dot{\xi}-3 \xi)=- Q(g,R)_{233424};$$
$$Q(g,S)_{1313}=\frac{t (t \ddot{\xi}-2 \dot{\xi})}{t-2 \xi}, \ \ Q(g,S)_{1414}=\frac{t \sin^2\theta (t \ddot{\xi}-2 \dot{\xi})}{t-2 \xi},$$
$$Q(g,S)_{2323}=-\frac{(t-2 \xi) (t \ddot{\xi}-2 \dot{\xi})}{t}, \ \ Q(g,S)_{2424}=-\frac{(t-2 \xi) \sin^2\theta (t \ddot{\xi}-2 \dot{\xi})}{t};$$
$$Q(g,C)_{122313}=-\frac{1}{2} t \ddot{\xi}+2 \dot{\xi}-\frac{3 \xi}{t} = -Q(g,C)_{121323},$$
$$Q(g,C)_{143413}=\frac{t^2 \sin^2\theta \left(t (t \ddot{\xi}-4 \dot{\xi})+6 \xi\right)}{2 (t-2 \xi)} = -Q(g,C)_{133414},$$
$$Q(g,C)_{121424}=\frac{\sin^2\theta \left(t (t \ddot{\xi}-4 \dot{\xi})+6\xi\right)}{2 t} = -Q(g,C)_{122414},$$
$$Q(g,C)_{233424}=\frac{1}{2} (t-2 \xi) \sin^2\theta \left(t (t \ddot{\xi}-4 \dot{\xi})+6 \xi\right) = -Q(g,C)_{243423};$$
$$Q(g,P)_{131223} = Q(g,P)_{232113} = Q(g,P)_{121323} = - Q(g,P)_{122313}=t \ddot{\xi}-3 \dot{\xi}+\frac{3 \xi}{t},$$
$$Q(g,P)_{123213}= - Q(g,P)_{132123} = - Q(g,P)_{231213} = - Q(g,P)_{123123} =\frac{2}{3} t \ddot{\xi}-\frac{7 \dot{\xi}}{3}+\frac{3 \xi}{t},$$
$$Q(g,P)_{131113}=-\frac{t^2 (t \ddot{\xi}-2 \dot{\xi})}{3 (t-2 \xi)^2}, \ \ Q(g,P)_{133313}=\frac{t^3 (t \ddot{\xi}-2 \dot{\xi})}{3 (t-2 \xi)},$$
$$Q(g,P)_{143413} = Q(g,P)_{134314} = Q(g,P)_{341413} = - Q(g,P)_{341314} =\frac{t^2 \sin^2 \theta \left(t (t \ddot{\xi}-5 \dot{\xi})+9 \xi\right)}{3 (t-2 \xi)},$$
$$Q(g,P)_{144313} = Q(g,P)_{133414} = Q(g,P)_{344113} = - Q(g,P)_{343114} =\frac{t^2 \sin^2 \theta (t \dot{\xi}-3 \xi)}{t-2 \xi},$$
$$Q(g,P)_{141224} = Q(g,P)_{242114} = Q(g,P)_{121424} = - Q(g,P)_{122414}=\frac{\sin^2 \theta \left(t (t \ddot{\xi}-3 \dot{\xi})+3 \xi\right)}{t},$$
$$Q(g,P)_{142124} = Q(g,P)_{124124} = Q(g,P)_{241214} = - Q(g,P)_{124214}=-\frac{\sin^2 \theta \left(t (2 t \ddot{\xi}-7 \dot{\xi})+9 \xi\right)}{3 t},$$
$$Q(g,P)_{141114}=-\frac{t^2 \sin^2 \theta (t \ddot{\xi}-2 \dot{\xi})}{3 (t-2 \xi)^2}, \ \ Q(g,P)_{144414}=\frac{t^3 \sin^4 \theta (t \ddot{\xi}-2 \dot{\xi})}{3 (t-2 \xi)},$$
$$Q(g,P)_{232223}=-\frac{(t-2 \xi)^2 (t \ddot{\xi}-2 \dot{\xi})}{3 t^2}, \ \ Q(g,P)_{233323}=-\frac{1}{3} t (t-2 \xi) (t \ddot{\xi}-2 \dot{\xi}),$$
$$Q(g,P)_{243423} = Q(g,P)_{342423} = Q(g,P)_{234324}$$
$$ = - Q(g,P)_{342324} =-\frac{1}{3} (t-2 \xi) \sin^2 \theta \left(t (t \ddot{\xi}-5 \dot{\xi})+9 \xi\right),$$
$$Q(g,P)_{244323} = Q(g,P)_{344223} = Q(g,P)_{233424} = - Q(g,P)_{343224} =-(t-2 \xi) \sin^2 \theta (t \dot{\xi}-3 \xi),$$
$$Q(g,P)_{242224}=-\frac{(t-2 \xi)^2 \sin^2 \theta (t \ddot{\xi}-2 \dot{\xi})}{3 t^2}, \ \ Q(g,P)_{244424}=-\frac{1}{3} t (t-2 \xi) \sin^4 \theta (t \ddot{\xi}-2 \dot{\xi}).$$
\newline
\noindent
{\bf{Part II.}} For the tensors $C \cdot R$, $C \cdot S$, $C \cdot C$ and $C \cdot P$ we have the following relations:
$$(C \cdot R)_{121323}=\frac{\left(t (t \ddot{\xi}-4 \dot{\xi})+6 \xi\right) \left(t (t \ddot{\xi}-3 \dot{\xi})+3 \xi\right)}{6 t^4} = -(C \cdot R)_{122313},$$
$$(C \cdot R)_{133414}=\frac{\sin^2\theta (t \dot{\xi}-3 \xi) \left(t (t \ddot{\xi}-4 \dot{\xi})+6 \xi\right)}{6 t (t-2 \xi)} = -(C \cdot R)_{143413},$$
$$(C \cdot R)_{121424}=\frac{\sin^2\theta \left(t (t \ddot{\xi}-4 \dot{\xi})+6 \xi\right) \left(t (t \ddot{\xi}-3 \dot{\xi})+3 \xi\right)}{6 t^4} = -(C \cdot R)_{122414},$$
$$(C \cdot R)_{243423}=\frac{(t-2 \xi) \sin^2\theta (t \dot{\xi}-3\xi) \left(t (t \ddot{\xi}-4 \dot{\xi})+6 \xi\right)}{6 t^3} = -(C \cdot R)_{233424};$$
$$(C \cdot S)_{1313}=\frac{(t \ddot{\xi}-2 \dot{\xi}) \left(t(t \ddot{\xi}-4 \dot{\xi})+6\xi\right)}{6 t^2 (t-2 \xi)},$$
$$(C \cdot S)_{1414}=\frac{\sin^2\theta (t\ddot{\xi}-2 \dot{\xi}) \left(t (t \ddot{\xi}-4 \dot{\xi})+6\xi\right)}{6 t^2 (t-2 \xi)},$$
$$(C \cdot S)_{2323}=-\frac{(t-2 \xi) (t\ddot{\xi}-2 \dot{\xi}) \left(t (t \ddot{\xi}-4 \dot{\xi})+6 \xi\right)}{6 t^4},$$
$$(C \cdot S)_{2424}=-\frac{(t-2 \xi) \sin^2\theta (t \ddot{\xi}-2 \dot{\xi}) \left(t (t \ddot{\xi}-4 \dot{\xi})+6 \xi\right)}{6 t^4};$$
$$(C \cdot C)_{121323}=\frac{\left(t (t \ddot{\xi}-4 \dot{\xi})+6 \xi\right)^2}{12 t^4}= - (C \cdot C)_{122313},$$
$$(C \cdot C)_{143413}=\frac{\sin^2\theta \left(t (t \ddot{\xi}-4 \dot{\xi})+6 \xi\right)^2}{12 t (t-2 \xi)} = -(C \cdot C)_{133414},$$
$$(C \cdot C)_{121424}=\frac{\sin^2\theta \left(t (t \ddot{\xi}-4 \dot{\xi})+6 \xi\right)^2}{12 t^4} = -(C \cdot C)_{122414},$$
$$(C \cdot C)_{233424}=\frac{(t-2 \xi) \sin^2\theta \left(t (t \ddot{\xi}-4 \dot{\xi})+6 \xi\right)^2}{12 t^3}= -(C \cdot C)_{243423};$$
$$(C \cdot P)_{131223} = (C \cdot P)_{231213} = (C \cdot P)_{121323}$$
$$ = - (C \cdot P)_{122313}=\frac{\left(t (t \ddot{\xi}-4 \dot{\xi})+6 \xi\right) \left(t (t \ddot{\xi}-3 \dot{\xi})+3 \xi\right)}{6 t^4},$$
$$(C \cdot P)_{132123} = (C \cdot P)_{231213} = (C \cdot P)_{123123}$$
$$ = - (C \cdot P)_{123213}= - \frac{\left(t (t \ddot{\xi}-4 \dot{\xi})+6 \xi\right) \left(t (2 t \ddot{\xi}-7 \dot{\xi})+9 \xi\right)}{18 t^4},$$
$$(C \cdot P)_{131113}=-\frac{(t \ddot{\xi}-2 \dot{\xi}) \left(t (t \ddot{\xi}-4 \dot{\xi})+6 \xi\right)}{18 t (t-2 \xi)^2}, \ \ (C \cdot P)_{133313}=\frac{(t \ddot{\xi}-2 \dot{\xi}) \left(t (t \ddot{\xi}-4 \dot{\xi})+6 \xi\right)}{18 (t-2 \xi)},$$
$$(C \cdot P)_{143413} = (C \cdot P)_{341413} = (C \cdot P)_{134314}$$
$$ = - (C \cdot P)_{341314} =\frac{\sin^2\theta \left(t (t \ddot{\xi}-5 \dot{\xi})+9 \xi\right) \left(t (t \ddot{\xi}-4 \dot{\xi})+6 \xi\right)}{18 t(t-2 \xi)},$$
$$(C \cdot P)_{144313} = (C \cdot P)_{344113} = (C \cdot P)_{133414} = - (C \cdot P)_{343114} =\frac{\sin^2\theta (t \dot{\xi}-3 \xi) \left(t (t \ddot{\xi}-4 \dot{\xi})+6 \xi\right)}{6 t (t-2 \xi)},$$
$$(C \cdot P)_{141224} = (C \cdot P)_{242114} = (C \cdot P)_{121424}$$
$$ = - (C \cdot P)_{122414}=\frac{\sin^2\theta \left(t (t \ddot{\xi}-4 \dot{\xi})+6 \xi\right) \left(t (t \ddot{\xi}-3 \dot{\xi})+3 \xi\right)}{6 t^4},$$
$$(C \cdot P)_{142124} = (C \cdot P)_{241214} = (C \cdot P)_{124124}$$
$$ = - (C \cdot P)_{124214}= -\frac{\sin^2\theta \left(t (t \ddot{\xi}-4 \dot{\xi})+6 \xi\right) \left(t (2 t \ddot{\xi}-7 \dot{\xi})+9 \xi\right)}{18 t^4},$$
$$(C \cdot P)_{141114}=-\frac{\sin^2\theta (t \ddot{\xi}-2 \dot{\xi}) \left(t (t \ddot{\xi}-4 \dot{\xi})+6 \xi\right)}{18 t (t-2 \xi)^2},$$
$$(C \cdot P)_{144414} = \frac{\sin^4\theta (t \ddot{\xi}-2 \dot{\xi}) \left(t (t \ddot{\xi}-4 \dot{\xi})+6 \xi\right)}{18 (t-2 \xi)},$$
$$(C \cdot P)_{232223}=-\frac{(t-2 \xi)^2 (t \ddot{\xi}-2 \dot{\xi}) \left(t (t \ddot{\xi}-4 \dot{\xi})+6 \xi\right)}{18 t^5},$$
$$(C \cdot P)_{233323}=-\frac{(t-2 \xi) (t \ddot{\xi}-2 \dot{\xi}) \left(t (t \ddot{\xi}-4 \dot{\xi})+6 \xi\right)}{18 t^2},$$
$$(C \cdot P)_{243423} = (C \cdot P)_{342423} = (C \cdot P)_{234324}$$
$$ = - (C \cdot P)_{342324} =-\frac{(t-2 \xi) \sin^2\theta \left(t (t \ddot{\xi}-5 \dot{\xi})+9 \xi\right) \left(t (t \ddot{\xi}-4 \dot{\xi})+6 \xi\right)}{18 t^3},$$
$$(C \cdot P)_{244323} = (C \cdot P)_{344223} = (C \cdot P)_{233424}$$
$$ = - (C \cdot P)_{343224} =-\frac{(t-2 \xi) \sin^2\theta (t \dot{\xi}-3 \xi) \left(t (t \ddot{\xi}-4 \dot{\xi})+6 \xi\right)}{6 t^3},$$
$$(C \cdot P)_{242224}=-\frac{(t-2 \xi)^2 \sin^2\theta (t \ddot{\xi}-2 \dot{\xi}) \left(t (t \ddot{\xi}-4 \dot{\xi})+6 \xi\right)}{18 t^5},$$
$$(C \cdot P)_{244424}=-\frac{(t-2 \xi) \sin^4\theta (t \ddot{\xi}-2 \dot{\xi}) \left(t (t \ddot{\xi}-4 \dot{\xi})+6 \xi\right)}{18 t^2}.$$

For the tensors $W \cdot R$, $W \cdot S$, $W \cdot C$ and $W \cdot P$ we have the following relations:
$$(W \cdot R)_{121323}=\frac{\left(t (t \ddot{\xi}-4 \dot{\xi})+6 \xi\right) \left(t (t \ddot{\xi}-3\dot{\xi})+3 \xi\right)}{6 t^4} = -(W \cdot R)_{122313},$$
$$(W \cdot R)_{133414}=\frac{\sin^2\theta (t \dot{\xi}-3 \xi)\left(t (t \ddot{\xi}-4 \dot{\xi})+6 \xi\right)}{6 t (t-2 \xi)}= -(W \cdot R)_{143413},$$
$$(W \cdot R)_{121424}=\frac{\sin^2\theta \left(t (t \ddot{\xi}-4 \dot{\xi})+6 \xi\right) \left(t(t \ddot{\xi}-3 \dot{\xi})+3 \xi\right)}{6 t^4} = -(W \cdot R)_{122414},$$
$$(W \cdot R)_{243423}=\frac{(t-2 \xi) \sin^2\theta (t \dot{\xi}-3 \xi) \left(t (t \ddot{\xi}-4 \dot{\xi})+6 \xi\right)}{6 t^3}= - (W \cdot R)_{233424};$$
$$(W \cdot S)_{1313}=\frac{(t \ddot{\xi}-2 \dot{\xi}) \left(t (t \ddot{\xi}-4 \dot{\xi})+6 \xi\right)}{6 t^2 (t-2 \xi)},$$
$$(W \cdot S)_{1414}=\frac{\sin^2\theta (t \ddot{\xi}-2 \dot{\xi}) \left(t(t \ddot{\xi}-4 \dot{\xi})+6 \xi\right)}{6 t^2 (t-2 \xi)},$$
$$(W \cdot S)_{2323}=-\frac{(t-2 \xi) (t \ddot{\xi}-2 \dot{\xi}) \left(t (t\ddot{\xi}-4 \dot{\xi})+6 \xi\right)}{6 t^4},$$
$$(W \cdot S)_{2424}=-\frac{(t-2 \xi) \sin^2\theta (t \ddot{\xi}-2 \dot{\xi}) \left(t (t \ddot{\xi}-4 \dot{\xi})+6 \xi\right)}{6 t^4};$$
$$(W \cdot C)_{121323}=\frac{\left(t (t \ddot{\xi}-4 \dot{\xi})+6 \xi\right)^2}{12 t^4}= -(W \cdot C)_{122313},$$
$$(W \cdot C)_{143413}=\frac{\sin^2\theta \left(t (t \ddot{\xi}-4 \dot{\xi})+6 \xi\right)^2}{12 t (t-2 \xi)} =-(W \cdot C)_{133414},$$
$$(W \cdot C)_{121424}=\frac{\sin^2\theta \left(t (t \ddot{\xi}-4 \dot{\xi})+6 \xi\right)^2}{12 t^4}= -(W \cdot C)_{122414},$$
$$(W \cdot C)_{233424}=\frac{(t-2 \xi) \sin^2\theta \left(t (t \ddot{\xi}-4 \dot{\xi})+6 \xi\right)^2}{12 t^3} =-(W \cdot C)_{243423};$$
$$(W \cdot P)_{121323} = - (W \cdot P)_{122313} =\frac{\left(t (t \ddot{\xi}-4 \dot{\xi})+6 \xi\right) \left(t (t \ddot{\xi}-3 \dot{\xi})+3 \xi\right)}{6 t^4},$$
$$(W \cdot P)_{123213} = - (W \cdot P)_{123123} = - (W \cdot P)_{132123} =\frac{\left(t (t \ddot{\xi}-4 \dot{\xi})+6 \xi\right) \left(t (2 t \ddot{\xi}-7 \dot{\xi})+9 \xi\right)}{18 t^4},$$
$$(W \cdot P)_{133313} = \frac{1}{\sin^4 \theta}(W \cdot P)_{144414} =\frac{(t \ddot{\xi}-2 \dot{\xi}) \left(t (t \ddot{\xi}-4 \dot{\xi})+6 \xi\right)}{18 (t-2 \xi)},$$
$$(W \cdot P)_{143413} = (W \cdot P)_{134314} =\frac{\sin^2 \theta \left(t (t \ddot{\xi}-5 \dot{\xi})+9 \xi\right) \left(t (t \ddot{\xi}-4 \dot{\xi})+6 \xi\right)}{18 t (t-2 \xi)},$$
$$(W \cdot P)_{144313} = (W \cdot P)_{133414} = (W \cdot P)_{344113} =\frac{\sin^2 \theta (t \dot{\xi}-3 \xi) \left(t (t \ddot{\xi}-4 \dot{\xi})+6 \xi\right)}{6 t (t-2 \xi)},$$
$$(W \cdot P)_{122414} = - (W \cdot P)_{121424} =-\frac{\sin^2 \theta \left(t (t \ddot{\xi}-4 \dot{\xi})+6 \xi\right) \left(t (t \ddot{\xi}-3 \dot{\xi})+3 \xi\right)}{6 t^4},$$
$$(W \cdot P)_{124214} = - (W \cdot P)_{142124} = - (W \cdot P)_{124124} =\frac{\sin^2 \theta \left(t (t \ddot{\xi}-4 \dot{\xi})+6 \xi\right) \left(t (2 t \ddot{\xi}-7 \dot{\xi})+9 \xi\right)}{18 t^4},$$
$$(W \cdot P)_{233323} = \frac{1}{\sin^4 \theta}(W \cdot P)_{244424} =-\frac{(t-2 \xi) (t \ddot{\xi}-2 \dot{\xi}) \left(t (t \ddot{\xi}-4 \dot{\xi})+6 \xi\right)}{18 t^2},$$
$$(W \cdot P)_{243423} = (W \cdot P)_{234324} =-\frac{(t-2 \xi) \sin^2 \theta \left(t (t \ddot{\xi}-5 \dot{\xi})+9 \xi\right) \left(t (t \ddot{\xi}-4 \dot{\xi})+6 \xi\right)}{18 t^3},$$
$$(W \cdot P)_{244323} = (W \cdot P)_{344223} = (W \cdot P)_{233424} =-\frac{(t-2 \xi) \sin^2 \theta (t \dot{\xi}-3 \xi) \left(t (t \ddot{\xi}-4 \dot{\xi})+6 \xi\right)}{6 t^3}.$$

For the tensors 
$K \cdot R$, $K \cdot S$, $K \cdot C$, $K \cdot W$, $K \cdot K$ and $K \cdot P$ we have the following relations:
$$(K \cdot R)_{121323}=\frac{(t^2 \ddot{\xi}+2 \xi) \left(t (t \ddot{\xi}-3 \dot{\xi})+3 \xi\right)}{2 t^4} = -(K \cdot R)_{122313},$$
$$(K \cdot R)_{133414}=\frac{\sin^2\theta (t \dot{\xi}-3 \xi) (t^2 \ddot{\xi}+2 \xi)}{2 t (t-2 \xi)} = -(K \cdot R)_{143413},$$
$$(K \cdot R)_{121424}=\frac{\sin^2\theta (t^2 \ddot{\xi}+2 \xi) \left(t (t \ddot{\xi}-3 \dot{\xi})+3 \xi\right)}{2 t^4} = -(K \cdot R)_{122414},$$
$$(K \cdot R)_{243423}=\frac{(t-2 \xi) \sin^2\theta (t \dot{\xi}-3 \xi) (t^2 \ddot{\xi}+2 \xi)}{2 t^3} = -(K \cdot R)_{233424};$$
$$(K \cdot S)_{1313}=\frac{(t^2 \ddot{\xi}+2 \xi) (t\ddot{\xi}-2 \dot{\xi})}{2 t^2 (t-2 \xi)}, \ \ (K \cdot S)_{1414}=\frac{\sin^2\theta (t^2 \ddot{\xi}+2 \xi) (t \ddot{\xi}-2 \dot{\xi})}{2 t^2 (t-2 \xi)},$$
$$(K \cdot S)_{2323}=-\frac{(t-2 \xi) (t^2 \ddot{\xi}+2 \xi) (t \ddot{\xi}-2 \dot{\xi})}{2 t^4}, \ \ (K \cdot S)_{2424}=-\frac{(t-2 \xi) \sin^2\theta (t^2 \ddot{\xi}+2 \xi) (t \ddot{\xi}-2 \dot{\xi})}{2 t^4};$$
$$(K \cdot C)_{121323}=\frac{(t^2 \ddot{\xi}+2 \xi) \left(t (t \ddot{\xi}-4 \dot{\xi})+6 \xi\right)}{4 t^4} = -(K \cdot C)_{122313},$$
$$(K \cdot C)_{143413}=\frac{\sin^2\theta (t^2 \ddot{\xi}+2 \xi) \left(t (t \ddot{\xi}-4 \dot{\xi})+6 \xi\right)}{4 t (t-2 \xi)} = -(K \cdot C)_{133414},$$
$$(K \cdot C)_{121424}=\frac{\sin^2\theta (t^2 \ddot{\xi}+2 \xi) \left(t (t \ddot{\xi}-4 \dot{\xi})+6 \xi\right)}{4 t^4} = -(K \cdot C)_{122414},$$
$$(K \cdot C)_{233424}=\frac{(t-2 \xi) \sin^2\theta (t^2 \ddot{\xi}+2 \xi) \left(t (t \ddot{\xi}-4 \dot{\xi})+6 \xi\right)}{4 t^3} = -(K \cdot C)_{243423};$$
$$(K \cdot P)_{122313} = -(K \cdot P)_{121323} = \frac{1}{\sin^2 \theta}(K \cdot P)_{122414}$$
$$ = - \frac{1}{\sin^2 \theta}(K \cdot P)_{121424} = -\frac{(t^2 \ddot{\xi}+2 \xi) \left(t (t \ddot{\xi}-3 \dot{\xi})+3 \xi\right)}{2 t^4},$$
$$(K \cdot P)_{123213} = - (K \cdot P)_{132123} = - (K \cdot P)_{123123} =\frac{(t^2 \ddot{\xi}+2 \xi) \left(t (2 t \ddot{\xi}-7 \dot{\xi})+9 \xi\right)}{6 t^4},$$
$$(K \cdot P)_{133313} = \frac{1}{\sin^4 \theta}(K \cdot P)_{144414} =\frac{(t^2 \ddot{\xi}+2 \xi) (t \ddot{\xi}-2 \dot{\xi})}{6 (t-2 \xi)},$$
$$(K \cdot P)_{143413} = (K \cdot P)_{134314} =\frac{\sin^2\theta (t^2 \ddot{\xi}+2 \xi) \left(t(t\ddot{\xi}-5 \dot{\xi})+9 \xi\right)}{6 t (t-2 \xi)},$$
$$(K \cdot P)_{144313} = (K \cdot P)_{133414} = (K \cdot P)_{344113} =\frac{\sin^2\theta (t \dot{\xi}-3 \xi) (t^2 \ddot{\xi}+2 \xi)}{2 t (t-2 \xi)},$$
$$(K \cdot P)_{233323} = \frac{1}{\sin^4 \theta}(K \cdot P)_{244424} =-\frac{(t-2 \xi) (t^2 \ddot{\xi}+2 \xi) (t \ddot{\xi}-2 \dot{\xi})}{6 t^2},$$
$$(K \cdot P)_{244323} = (K \cdot P)_{344223} = (K \cdot P)_{233424} =-\frac{(t-2 \xi) \sin^2\theta (t \dot{\xi}-3 \xi) (t^2 \ddot{\xi}+2 \xi)}{2 t^3},$$
$$(K \cdot P)_{124124} = (K \cdot P)_{142124} = -(K \cdot P)_{124214} =-\frac{\sin^2\theta (t^2 \ddot{\xi}+2 \xi) \left(t (2 t \ddot{\xi}-7 \dot{\xi})+9 \xi\right)}{6 t^4},$$
$$(K \cdot P)_{234324} = (K \cdot P)_{243423} =-\frac{(t-2 \xi) \sin^2\theta (t^2 \ddot{\xi}+2 \xi) \left(t (t \ddot{\xi}-5 \dot{\xi})+9 \xi \right)}{6 t^3}.$$

For the tensors
$P \cdot R$, $P \cdot S$, $P \cdot C$, $P \cdot W$, $P \cdot K$ and $P \cdot P$ we have the following relations:
$$(P \cdot R)_{121323}=\frac{t^2 \dot{\xi} \left(5 \dot{\xi}-2 t \ddot{\xi}\right)+2 t \xi \left(2 t \ddot{\xi}-7 \dot{\xi}\right)+9 \xi^2}{3 t^4}=(P \cdot R)_{122331},$$
$$(P \cdot R)_{121424}=\frac{\sin^2\theta \left(t^2 \dot{\xi} \left(5 \dot{\xi}-2 t \ddot{\xi}\right)+2 t \xi \left(2 t \ddot{\xi}-7 \dot{\xi}\right)+9 \xi^2\right)}{3 t^4}=(P \cdot R)_{124214},$$
$$(P \cdot R)_{143413}=\frac{\sin^2\theta \left(t^2 \dot{\xi}^2+2 t \xi \left(t \ddot{\xi}-5 \dot{\xi}\right)+9 \xi^2\right)}{3 t (t-2 \xi)}=(P \cdot R)_{134314},$$
$$(P \cdot R)_{233424}=\frac{(t-2 \xi) \sin^2\theta \left(t^2 \dot{\xi}^2+2 t \xi \left(t \ddot{\xi}-5 \dot{\xi}\right)+9 \xi^2\right)}{3 t^3}=(P \cdot R)_{244323};$$
$$(P \cdot S)_{1313}=\frac{(\xi-t \dot{\xi}) (t \ddot{\xi}-2 \dot{\xi})}{t^2 (t-2 \xi)} = -(P \cdot S)_{1331},$$
$$(P \cdot S)_{1414}=\frac{\sin^2\theta (\xi-t \dot{\xi}) (t \ddot{\xi}-2 \dot{\xi})}{t^2 (t-2 \xi)} = -(P \cdot S)_{1441},$$
$$(P \cdot S)_{2323}=\frac{(t-2 \xi) (t \dot{\xi}-\xi) (t \ddot{\xi}-2 \dot{\xi})}{t^4} = - (P \cdot S)_{2332},$$
$$(P \cdot S)_{2424}=\frac{(t-2 \xi) \sin^2\theta (t \dot{\xi}-\xi) (t \ddot{\xi}-2 \dot{\xi})}{t^4}= - (P \cdot S)_{2442};$$
$$(P \cdot C)_{121323}=\frac{\left(t \left(t \ddot{\xi}-5 \dot{\xi}\right)+9 \xi\right) \left(t \left(t \ddot{\xi}-4 \dot{\xi}\right)+6 \xi\right)}{18 t^4}=(P \cdot C)_{123213},$$
$$(P \cdot C)_{121424}=\frac{\sin^2\theta \left(t \left(t \ddot{\xi}-5 \dot{\xi}\right)+9 \xi\right) \left(t \left(t \ddot{\xi}-4 \dot{\xi}\right)+6 \xi\right)}{18 t^4}=(P \cdot C)_{124214},$$
$$(P \cdot C)_{134314}=\frac{\sin^2\theta \left(t \left(t \ddot{\xi}-4 \dot{\xi}\right)+6 \xi\right) \left(t \left(2 t \ddot{\xi}-7 \dot{\xi}\right)+9 \xi\right)}{18 t (t-2 \xi)}=(P \cdot C)_{143413},$$
$$(P \cdot C)_{233424}=\frac{(t-2 \xi) \sin^2\theta \left(t \left(t \ddot{\xi}-4 \dot{\xi}\right)+6 \xi\right) \left(t \left(2 t \ddot{\xi}-7 \dot{\xi}\right)+9 \xi\right)}{18 t^3}=(P \cdot C)_{244323};$$
$$(P \cdot P)_{121323}=\frac{\left(t \left(t \ddot{\xi}-5 \dot{\xi}\right)+9 \xi\right) \left(t \left(t \ddot{\xi}-3 \dot{\xi}\right)+3 \xi\right)}{9 t^4}=(P \cdot P)_{232113},$$
$$(P \cdot P)_{121424}=\frac{\sin^2\theta \left(t \left(t \ddot{\xi}-5 \dot{\xi}\right)+9 \xi\right) \left(t \left(t \ddot{\xi}-3 \dot{\xi}\right)+3 \xi\right)}{9 t^4}=(P \cdot P)_{242114},$$
$$(P \cdot P)_{123123}=\frac{\left(t \dot{\xi}-3 \xi\right) \left(t \left(t \ddot{\xi}-3 \dot{\xi}\right)+3 \xi\right)}{3 t^4}=(P \cdot P)_{231213},$$
$$(P \cdot P)_{124124}=\frac{\sin^2\theta \left(t \dot{\xi}-3 \xi\right) \left(t \left(t \ddot{\xi}-3 \dot{\xi}\right)+3 \xi\right)}{3 t^4}=(P \cdot P)_{241214},$$
$$(P \cdot P)_{131113}=-\frac{\left(t \ddot{\xi}-2 \dot{\xi}\right) \left(t \left(t \ddot{\xi}-3 \dot{\xi}\right)+3 \xi\right)}{9 t (t-2 \xi)^2}=-(P \cdot P)_{131131},$$
$$(P \cdot P)_{133331}=\frac{\left(t \dot{\xi}-3 \xi\right) \left(t \ddot{\xi}-2 \dot{\xi}\right)}{9 (t-2 \xi)}=-(P \cdot P)_{133313},$$
$$(P \cdot P)_{133414}=\frac{\sin^2\theta \left(t \dot{\xi}-3 \xi\right) \left(t \left(2 t \ddot{\xi}-7 \dot{\xi}\right)+9 \xi\right)}{9 t (t-2 \xi)}=(P \cdot P)_{344113},$$
$$(P \cdot P)_{134341}=\frac{\sin^2\theta \left(t \dot{\xi}-3 \xi\right) \left(t \left(t \ddot{\xi}-3 \dot{\xi}\right)+3 \xi\right)}{3 t (t-2\xi)}=(P \cdot P)_{341431},$$
$$(P \cdot P)_{141141}=\frac{\sin^2\theta \left(t \ddot{\xi}-2 \dot{\xi}\right) \left(t \left(t \ddot{\xi}-3 \dot{\xi}\right)+3 \xi\right)}{9 t (t-2 \xi)^2}=-(P \cdot P)_{141114},$$
$$(P \cdot P)_{144441}=\frac{\sin^4\theta \left(t \dot{\xi}-3 \xi\right) \left(t \ddot{\xi}-2 \dot{\xi}\right)}{9 (t-2 \xi)}=-(P \cdot P)_{144414},$$
$$(P \cdot P)_{232232}=\frac{(t-2 \xi)^2 \left(t \ddot{\xi}-2 \dot{\xi}\right) \left(t \left(t \ddot{\xi}-3 \dot{\xi}\right)+3 \xi\right)}{9 t^5}=-(P \cdot P)_{232223},$$
$$(P \cdot P)_{233323}=\frac{(t-2 \xi) \left(t \dot{\xi}-3 \xi\right) \left(t \ddot{\xi}-2 \dot{\xi}\right)}{9 t^2}=-(P \cdot P)_{233332},$$
$$(P \cdot P)_{233442}=\frac{(t-2 \xi) \sin^2\theta \left(t \dot{\xi}-3 \xi\right) \left(t \left(2 t \ddot{\xi}-7 \dot{\xi}\right)+9 \xi\right)}{9 t^3}=(P \cdot P)_{344232},$$
$$(P \cdot P)_{234324}=\frac{(t-2 \xi) \sin^2\theta \left(t \dot{\xi}-3 \xi\right) \left(t \left(t \ddot{\xi}-3 \dot{\xi}\right)+3 \xi\right)}{3 t^3}=(P \cdot P)_{342423},$$
$$(P \cdot P)_{242242}=\frac{(t-2 \xi)^2 \sin^2\theta \left(t \ddot{\xi}-2 \dot{\xi}\right) \left(t \left(t \ddot{\xi}-3 \dot{\xi}\right)+3 \xi\right)}{9 t^5}=-(P \cdot P)_{242224},$$
$$(P \cdot P)_{244424}=\frac{(t-2 \xi) \sin^4\theta \left(t \dot{\xi}-3 \xi\right) \left(t \ddot{\xi}-2 \dot{\xi}\right)}{9 t^2}=-(P \cdot P)_{244442}.$$

For the tensors $Q(S, R)$, $Q(S, C)$, $Q(S, W)$, $Q(S, K)$ and $Q(S, P)$ we have the following relations:
$$Q(S,R)_{122313}=\frac{4 \dot{\xi}(\xi-t \dot{\xi})+t (t \dot{\xi}+\xi) \ddot{\xi}}{t^3} = - Q(S,R)_{121323},$$
$$Q(S,R)_{133414}=-\frac{2 \sin^2\theta \left(\xi (t \ddot{\xi}+\dot{\xi})-t \dot{\xi}^2\right)}{t-2 \xi} = -Q(S,R)_{143413},$$
$$Q(S,R)_{122414}=\frac{\sin^2\theta \left(4 \dot{\xi} (\xi-t \dot{\xi})+t (t \dot{\xi}+\xi) \ddot{\xi}\right)}{t^3} = -Q(S,R)_{121424},$$
$$Q(S,R)_{243423}=\frac{2 (t-2 \xi) \sin^2\theta \left(\xi (t \ddot{\xi}+\dot{\xi})-t \dot{\xi}^2\right)}{t^2} = - Q(S,R)_{233424};$$
$$Q(S, C)_{121332}=\frac{\left(t \ddot{\xi}+4 \dot{\xi}\right) \left(t \left(t \ddot{\xi}-4 \dot{\xi}\right)+6 \xi\right)}{6 t^3}
=Q(S, C)_{122313},$$
$$Q(S, C)_{121442}=\frac{\sin^2\theta \left(t \ddot{\xi}+4 \dot{\xi}\right) \left(t \left(t \ddot{\xi}-4 \dot{\xi}\right)+6 \xi\right)}{6 t^3}=Q(S, C)_{122414},$$
$$Q(S, C)_{133414}=\frac{\sin^2\theta \left(t \ddot{\xi}+\dot{\xi}\right) \left(t \left(t \ddot{\xi}-4 \dot{\xi}\right)+6 \xi\right)}{3 (t-2 \xi)}=Q(S, C)_{144313},$$
$$Q(S, C)_{234324}=\frac{(t-2 \xi) \sin^2\theta \left(t \ddot{\xi}+\dot{\xi}\right) \left(t \left(t \ddot{\xi}-4 \dot{\xi}\right)+6 \xi\right)}{3 t^2}=Q(S, C)_{243423};$$
$$Q(S, W)_{121332}=\frac{t^3 \ddot{\xi}^2+4 \dot{\xi} \left(6 \xi-7 t \dot{\xi}\right)+6 t \left(t \dot{\xi}+\xi\right) \ddot{\xi}}{6 t^3}=Q(S, W)_{122313},$$
$$Q(S, W)_{121442}=\frac{\sin^2\theta \left(t^3 \ddot{\xi}^2+4 \dot{\xi} \left(6 \xi-7 t \dot{\xi}\right)+6 t \left(t \dot{\xi}+\xi\right) \ddot{\xi}\right)}{6 t^3}=Q(S, W)_{122414},$$
$$Q(S, W)_{133414}=\frac{\sin^2\theta \left(12 \xi \left(t \ddot{\xi}+\dot{\xi}\right)-t \left(t^2 \ddot{\xi}^2+8 \dot{\xi}^2\right)\right)}{6 (t-2 \xi)}=Q(S, W)_{143431},$$
$$Q(S, W)_{133441}=\frac{\sin^2\theta \left(t^3 \ddot{\xi}^2+8 t \dot{\xi}^2-12 \xi \left(t \ddot{\xi}+\dot{\xi}\right)\right)}{6 (t-2 \xi)}=Q(S, W)_{143413},$$
$$Q(S, W)_{233424}=\frac{(t-2 \xi) \sin^2\theta \left(t^3 \ddot{\xi}^2+8 t \dot{\xi}^2-12 \xi \left(t \ddot{\xi}+\dot{\xi}\right)\right)}{6 t^2}=Q(S, W)_{243432};$$
$$Q(S, K)_{121332}=\frac{t^3 \ddot{\xi}^2+2 t \xi \ddot{\xi}+8 \dot{\xi} \left(\xi-t \dot{\xi}\right)}{2 t^3}=Q(S, K)_{122313},$$
$$Q(S, K)_{121442}=\frac{\sin^2\theta \left(t^3 \ddot{\xi}^2+2 t \xi \ddot{\xi}+8 \dot{\xi} \left(\xi-t \dot{\xi}\right)\right)}{2 t^3}
=-Q(S, K)_{242114},$$
$$Q(S, K)_{133414}=\frac{\sin^2\theta \left(2\xi \dot{\xi}+t \left(2 \xi-t \dot{\xi}\right) \ddot{\xi}\right)}{t-2 \xi}
=Q(S, K)_{341431},$$
$$Q(S, K)_{133441}=\frac{\sin^2\theta \left(t \left(t \dot{\xi}-2 \xi\right) \ddot{\xi}-2 \xi \dot{\xi}\right)}{t-2 \xi}
=Q(S, K)_{344131},$$
$$Q(S, K)_{233424}=\frac{(t-2 \xi) \sin^2\theta \left(t \left(t \dot{\xi}-2 \xi\right) \ddot{\xi}-2 \xi \dot{\xi}\right)}{t^2}
=Q(S, K)_{342432};$$
$$Q(S, P)_{121332}=\frac{\left(t \ddot{\xi}+4 \dot{\xi}\right) \left(t \left(t \ddot{\xi}-3 \dot{\xi}\right)+3 \xi\right)}{3 t^3}=Q(S, P)_{122313}=-Q(S, P)_{131223},$$
$$Q(S, P)_{121442}=\frac{\sin^2\theta \left(t\ddot{\xi}+4 \dot{\xi}\right) \left(t \left(t \ddot{\xi}-3 \dot{\xi}\right)+3 \xi\right)}{3 t^3}=Q(S, P)_{122414}=-Q(S, P)_{141224},$$
$$Q(S, P)_{123123}=\frac{4 \dot{\xi} \left(\xi-t \dot{\xi}\right)+t \left(t \dot{\xi}+\xi\right) \ddot{\xi}}{t^3}=Q(S, P)_{132132}=Q(S, P)_{231231},$$
$$Q(S, P)_{124124}=\frac{\sin^2\theta \left(4 \dot{\xi} \left(\xi-t \dot{\xi}\right)+t \left(t \dot{\xi}+\xi\right) \ddot{\xi}\right)}{t^3}=Q(S, P)_{124214}=Q(S, P)_{142142},$$
$$Q(S, P)_{131113}=\frac{t \ddot{\xi} \left(t \ddot{\xi}-2 \dot{\xi}\right)}{3 (t-2 \xi)^2}=-Q(S, P)_{131131},$$
$$Q(S, P)_{133331}=\frac{2 t \dot{\xi} \left(t \ddot{\xi}-2 \dot{\xi}\right)}{3 (t-2 \xi)}=-Q(S, P)_{133313},$$
$$Q(S, P)_{133441}=\frac{2 \sin^2\theta \left(t \dot{\xi}-3 \xi\right) \left(t \ddot{\xi}+\dot{\xi}\right)}{3 (t-2 \xi)}=-Q(S, P)_{144313}=Q(S, P)_{343114},$$
$$Q(S, P)_{134341}=\frac{2 \sin^2\theta \left(\xi \left(t \ddot{\xi}+\dot{\xi}\right)-t \dot{\xi}^2\right)}{t-2 \xi}=-Q(S, P)_{134314}=-Q(S, P)_{341413},$$
$$Q(S, P)_{141114}=\frac{t \sin^2\theta \ddot{\xi} \left(t \ddot{\xi}-2 \dot{\xi}\right)}{3 (t-2 \xi)^2}=-Q(S, P)_{141141},$$
$$Q(S, P)_{144441}=\frac{2 t \sin^4\theta \dot{\xi} \left(t \ddot{\xi}-2 \dot{\xi}\right)}{3 (t-2 \xi)}=-Q(S, P)_{144414},$$
$$Q(S, P)_{232223}=\frac{(t-2 \xi)^2 \ddot{\xi} \left(t \ddot{\xi}-2 \dot{\xi}\right)}{3 t^3}=-Q(S, P)_{232232},$$
$$Q(S, P)_{233323}=\frac{2 (t-2 \xi) \dot{\xi} \left(t \ddot{\xi}-2 \dot{\xi}\right)}{3 t}=-Q(S, P)_{233332},$$
$$Q(S, P)_{233424}=\frac{2 (t-2 \xi) \sin^2\theta \left(t \dot{\xi}-3 \xi\right) \left(t \ddot{\xi}+\dot{\xi}\right)}{3 t^2}=Q(S, P)_{343242}=-Q(S, P)_{344232},$$
$$Q(S, P)_{234324}=\frac{2 (t-2 \xi) \sin^2\theta \left(\xi \left(t \ddot{\xi}+\dot{\xi}\right)-t \dot{\xi}^2\right)}{t^2}=-Q(S, P)_{342324}=-Q(S, P)_{243432},$$
$$Q(S, P)_{242224}=\frac{(t-2 \xi)^2 \sin^2\theta \ddot{\xi} \left(t \ddot{\xi}-2 \dot{\xi}\right)}{3 t^3}=-Q(S, P)_{242242},$$
$$Q(S, P)_{244424}=\frac{2 (t-2 \xi) \sin^4\theta \dot{\xi} \left(t \ddot{\xi}-2 \dot{\xi}\right)}{3 t}=-Q(S, P)_{244442}.$$

\noindent\textbf{Conclusion.} From the above results and discussion we conclude that Deszcz symmetric spaces 
are geometric models of the interior black hole spacetime and hence the defining conditions of Deszcz symmetric spaces are, physically, very stronger. 
However, for a specific value of $\xi$ the interior black hole spacetime turns into a semisymmetric spacetime and also a generalized Ricci pseudosymmetric spacetime. 

\noindent\textbf{Acknowledgement.} 
The work was carried out when the fourth named author visited Department of Mathematics of the Sardar Patel University 
as a visiting fellow under their UGC-SAP-DRS programme (F-510/5/DRS/2009 (SAP-II)). He also greatfully acknowledges 
the financial support of CSIR, New Delhi, India [Project F. No. 25(0171)/09/EMR-II].
The first named author
is supported by a grant of the Technische Universit\"at Berlin (Germany). 

\end{document}